\documentclass[english,12pt]{amsart}
\usepackage{amssymb,amsmath,amsthm}
\usepackage{babel,graphicx}

\newtheorem{theorem}{Theorem}[section]
\newtheorem{main}{Theorem}

\newtheorem{corollary}[theorem]{Corollary}
\newtheorem{lemma}[theorem]{Lemma}
\newtheorem{proposition}[theorem]{Proposition}
\theoremstyle{definition}
\newtheorem{definition}[theorem]{Definition}

\newtheorem{question}[theorem]{Question}

\theoremstyle{remark}
\newtheorem{remark}[theorem]{Remark}
\newtheorem{example}[theorem]{Example}
\numberwithin{equation}{section}

\newcommand{\fin}{\mbox{}\hfill $\square$ \\[0.2cm]}

\newenvironment{pf}{\par\noindent\textbf{Proof. }}
{\fin}
\newenvironment{pfof}[1]{\par\noindent
\textbf{Proof of #1. }}
{\fin}

\newcommand{\bbN}{\mathbb{N}}
\newcommand{\bbZ}{\mathbb{Z}}

\newcommand{\bbR}{\mathbb{R}}
\newcommand{\bbC}{\mathbb{C}}
\newcommand{\bfH}{\mathbf{H}}
\newcommand{\mcF}{\mathcal{F}}
\newcommand{\mcG}{\mathcal{G}}
\newcommand{\mcM}{\mathcal{M}}

\newcommand{\Laplace}{\Delta}
\newcommand{\grad}{\nabla}
\newcommand{\dist}{{\rm dist}}
\newcommand{\Diff}{{\rm Diff}}

\newcommand{\diam}{\mathop{\rm diam}}
\newcommand{\supp}{\mathop{\rm supp}}
\newcommand{\Expect}{{\mathbb E}}
\newcommand{\vol}{{\rm vol}}
\newcommand{\codim}{\mathop{\rm codim}}
\newcommand{\leb}{{\rm leb}}

\begin{document}

\title[Conformal dynamical systems]{Random conformal dynamical systems}
\author{Bertrand Deroin \& Victor Kleptsyn}%
\address{IHES - 35, Routes de Chartres - 91440 Bures sur Yvette - France}
\email{bderoin@ihes.fr}

\address{UMPA - CNRS UMR 5669 - ENS Lyon - 44, all\'ee d'Italie - 69007 Lyon - France}
\email{Victor.Kleptsyn@umpa.ens-lyon.fr}

\thanks{The first author acknowledges support from the Swiss National
Science Foundation.}
\thanks{The second author's work was partially supported by grants
RBFR 02-01-00482, RFBR 02-01-22002, CRDF RM1-2358-MO-02.}

\begin{abstract} We consider random dynamical systems such as groups of conformal transformations 
with a probability measure, or transversaly conformal foliations with a 
Laplace operator along the leaves, in which case we consider the  
holonomy pseudo-group. 
We prove that either there exists a measure invariant under all the elements of the group (or 
the pseudo-group), or almost surely a long composition of maps contracts exponentially a ball. 
We deduce some results about the unique ergodicity. 
\end{abstract}

\maketitle

\section{Introduction}

Given a probability measure on the group of homeomorphisms of a
manifold, one can study the asymptotic behaviour of large
composition of elements chosen randomly with respect to this
measure. From the 50's the questions of the behaviour of the random
walk on the group, composition of random matrices, equidistribution
of the orbits, Lyapunov exponents (if the action is by $C^1$
diffeomorphisms) etc. has been largely studied and understood. We
are unable to review here all the history of these problems and we
refer the interested reader to an excellent survey of Furman
\cite{Furman}. Nevertheless, we would like to mention here some
important works in the development of this theory:
Kakutani~\cite{Kakutani}, Furstenberg~\cite{Furstenberg},
Arnold-Krylov~\cite{Arnold-Krylov},
Furstenberg-Kesten~\cite{Furstenberg-Kesten}, Guivar'ch
\cite{Guivarc'h}.

\vspace{0.2cm}


More recently was developed the idea that the maps could be taken
from a pseudo-group, rather than a group. This has been introduced
in the paper of Garnett~\cite{Garnett} for the pseudo-group of a
foliation, and then studied by Ghys~\cite{Ghys2,Ghys3},
Kaimanovich~\cite{Kaimanovich}, Ledrappier~\cite{Ledrappier} and
Candel~\cite{Candel}. Following the lines of the ``Sullivan's
dictionary'', one can extend these ideas to other pseudo-groups,
like for instance the one generated by an endomorphism or a
correspondance.

\vspace{0.2cm}

In this work we study random compositions of the elements of a
pseudo-group acting conformally on a manifold. Our results concern
actions of a group by conformal transformations, conformal
correspondences or transversely conformal foliations (for one
dimensional dynamical systems, we suppose that the maps are of class
$C^1$); most of them were known in the case of a group acting
conformally on a compact manifold, when the conditional
probabilities do not depend on the point (see~\cite{Furman}).

\vspace{0.2cm}

Observe that the existence of a measure which is invariant by every
element of the pseudo-group is rather rare. In the symmetric case
(i.e. when the probabilities are symmetric) we prove the following
dichotomy. On a minimal subset, either is supported a probability
measure which is preserved by all the elements of the pseudo-group,
or the system has the property of exponential contraction: for any
point and almost every random composition of elements of the
pseudo-group there exists a neighborhood of the point which is
contracted exponentially.

\vspace{0.2cm}

We deduce some results about the equidistribution of the orbits of
the system along random compositions. In the case where the system
have the property of exponential contraction (even if the system is
not symmetric), we prove that the orbit of a point by almost every
random compositions is distributed with respect to a unique measure.
We also give examples of non symmetric systems for which the
exponential contraction property and the equidistribution property
are not satisfied.

\subsection{Presentation of the results for a foliation}
We begin by a survey on Garnett's theory~\cite{Garnett} (see also
\cite{Candel}).

\vspace{0.2cm}

Let $\mathcal F$ be a foliation of a compact manifold $M$, whose
leaves are of class $C^{\infty}$, and $g$ a Riemannian metric on the
leaves of~$\mathcal F$. In \cite{Garnett}, Garnett studies the
diffusion process along the leaves of~$\mathcal F$. Namely, the
metric $g$ induces the Laplace-Beltrami operator along the leaves,
which we denote $\Delta$; given a continuous function
$f_0:M\rightarrow \bbR$, one studies the heat equation along the
leaves of the foliation
\[  \frac{\partial f}{\partial t} = \Delta f\]
with initial condition $f(\cdot,t) = f_0$. As it is well-known
\cite{Chavel}, because the leaves are complete and of bounded
geometry, the solution to the heat equation is unique, defined for
all positive time, and is expressed by convolution of the initial
condition with the heat kernel $p(x,y;t)$. A fundamental Lemma due
to Garnett (see also Candel~\cite{Candel}), asserts that the
functions $f(\cdot,t)$ on $M$ are continuous, and that the diffusion
semi-group of operators $D^t$ defined for all $t\geq 0$ by
\[ D^t f(0,\cdot) = f(t,\cdot) \]
acts continuously on $C^0(M)$.

\vspace{0.2cm}

Associated to this diffusion semi-group, Garnett considers the
Brownian motion along the leaves of the foliation: this is a
Markovian process with continuous time, whose trajectories stay
every time in the same leaf, and whose transition probability
distributions are volume forms with leafwise density given by the heat
kernel. It is known that this process can be realized as a process
with continuous trajectories. For any point $x$ in $M$, let
$\Gamma_x$ be the set of continuous paths parametrized by
$[0,\infty)$, starting at $x$, and whose image is contained in the
leaf $\mcF_x$ passing through the point $x$. The space $\Gamma_x$ is
equipped with the uniform topology on compact subsets; there is a
probability Borel measure induced by the Brownian motion process,
expressing the probability that a trajectory occurs. This
probability measure is called the \textit{Wiener measure} and
denoted $W_x$.

\vspace{0.2cm}

We recall the definition of the \textit{holonomy pseudo-group}.
Because the manifold $M$ is compact, there is a
finite cover of $M$ by \textit{foliated box} $B_i \times T_i$, in which the
foliation $\mathcal F$ is the horizontal fibration. The change of coordinates
from $B_i\times T_i$ to $B_j\times T_j$ are of the form
\[  (x_i,t_i)\mapsto (x_j=x_j(x_i,t_i), t_j(t_i)).\]
The maps $t_j(t_i)$ generates a pseudo-group on the union $T=\cup_i T_i$, 
called the holonomy pseudo-group.
A measure on $T$ which is invariant by the holonomy pseudo-group is
called a \textit{transversely invariant measure}.
These measures have been
introduced by Schwartzman for flows~\cite{Schwartzman}, by Plante
and Ruelle-Sullivan for foliations~\cite{Plante,Ruelle-Sullivan} and
by Sullivan for other kind of dynamical systems~\cite{cycle}.

\vspace{0.2cm}

Now, consider a continuous path $\gamma$ contained in a leaf,
parametrized by a closed interval. It crosses successively the
foliation boxes $B_{i_1}\times T_{i_1},\ldots,B_{i_k}\times
T_{i_k}$. The composition of the associated change of transverse
coordinates is by definition the holonomy map $h_{\gamma}$
corresponding to $\gamma$. The following result describes the
asymptotic behaviour of the holonomy maps $h_{\gamma|_{[0,t]}}$ when
$t$ goes to infinity, for a generic Brownian path along the leaf
passing throw a point $x$, when the foliation is transversely
conformal.

\begin{theorem}[Main Theorem]\label{t:main theorem}
Let $\mathcal F$ be a transversely conformal foliation of class
$C^1$ of a compact manifold. Then either there exists a transversely
invariant measure. Or $\mathcal F$ has a finite number of minimal
sets $\mathcal M_1,\ldots,\mathcal M_k$ equipped with probability
measures $\mu_1,\ldots,\mu_k$, and there exists a real $\alpha >0$
such that:
\begin{itemize}
\item \textbf{Contraction.}
For every point $x$ in $M$ and almost every leafwise Brownian path
$\gamma$ starting at $x$, there is a neighborhood $T_{\gamma}$ of
$x$ in $T$ and a constant $C_{\gamma}>0$, such that for every $t>0$,
the holonomy map $h_{\gamma|_{[0,t]}}$ is defined on $T_{\gamma}$
and
\[   |h_{\gamma|_{[0,t]}}(T_{\gamma})| \leq C_{\gamma}\exp(-\alpha t).\]
\item \textbf{Distribution.}
For every point $x$ in $M$ and almost every leafwise Brownian path
$\gamma$ starting at $x$, the path $\gamma$ tends to one of the
$\mathcal M_j$ and is distributed with respect to $\mu_j$, in the
sense that
\[  \lim_{t\rightarrow \infty} \frac{1}{t} \gamma_* \leb_{[0,t]} = \mu_j,\]
where $\leb_{[0,t]}$ is the Lebesgue measure on the interval
$[0,t]$.
\item \textbf{Attraction.}
The probability $p_j(x)$ that a leafwise Brownian path starting at a
point $x$ of $M$ tends to $\mathcal M_j$ is a continuous leafwise
harmonic function.
\item \textbf{Diffusion.}
When $t$ goes to infinity, the diffusions $D^t f$ of a continuous
function $f:M\rightarrow\bbR$ converge uniformly to the function
$\sum_j c_j p_j$, where $c_j = \int f d\mu_j$. In particular, the
functions~$p_j$ form a base in the space of continuous leafwise
harmonic functions.
\end{itemize}
\end{theorem}

The existence of a transversely invariant measure for a transversely
conformal foliation is a very strong condition. An ergodic component
of such a measure is either supported on a compact leaf, or it is
diffuse. All the examples we know of a transversely conformal
foliation having a diffuse transversely invariant measure have a
transverse metric which is transversely invariant. In the case of
codimension one foliation of class $C^2$, this is an easy
consequence of Sacksteder Theorem. For higher transversely conformal
foliations this has been conjectured by Ghys~\cite{Ghys1}, and for
codimension 3 and higher with an additional restriction of
minimality this conjecture was proven by Tarquini~\cite{Tarquini}.

\vspace{0.2cm}

The contraction property for group actions on the circle was studied in \cite{Antonov},
\cite{Kaijser} and \cite{Kleptsyn-Nalski}.

\vspace{0.2cm} 

The distribution part of the theorem was proved by Garnett
\cite{Garnett} in the case of the stable foliation of the geodesic
flow on the unitary tangent bundle of a surface of constant negative
curvature, using the contraction property and the similarity, which
is straightforward in this case. This was also extended to the case
of a manifold of negative variable curvature of any dimension by
Ledrappier~\cite{Ledrappier} (also see \cite{Ledrappier2} and \cite{Yue}). 
Finally, Hamenst\"adt \cite{Hamenstadt1} has studied 
drifted Brownian motions on the stable foliation in negative curvature, and has obtained 
a sufficient condition for uniqueness of a harmonic measure, also giving an example 
(see \cite{Hamenstadt2}) of non uniqueness of the harmonic measure in the drifted case. 

\vspace{0.2cm}

The diffusion part of Theorem~\ref{t:main theorem} gives examples
of foliated riemannian manifolds that have non trivial continuous
leafwise harmonic functions. Such an example were constructed
in~\cite{Feres-Ghani}.

\subsection{Organization of the proof.}
In \cite{Garnett}, Garnett studies the ergodic properties of the
leafwise diffusion semi-group and of the leafwise Brownian motion.
She introduces the notion of \textit{harmonic measure}, which is a
probability measure invariant by the diffusion semi-group. By the
Kakutani fixed point Theorem, such a measure exists. The relation
with the Brownian motion goes as follows. Consider the space
$\Gamma$ of all the continuous paths contained in a leaf of
$\mathcal F$. There is a semi-group $\{\sigma_t\}_{t\geq 0}$ of
transformations of $\Gamma$ defined for every $t,s\geq 0$ and every
$\gamma \in \Gamma$ by
\[  \sigma_t(\gamma) (s) = \gamma(t+s).\]
If $\mu$ is a measure
on $M$, then one consider the probability measure $\overline{\mu}$ on $\Gamma$ which is defined
by
\[  \overline{\mu} (B) = \int _M W_x(B\cap \Gamma_x) d\mu(x),\]
for every Borel subset $B$ contained in $\Gamma$ (recall that for
every $x\in M$, $\Gamma_x$ is the set of pathes starting at $x$ and
contained on the leaf through $x$, and $W_x$ is the Wiener measure
on $\Gamma_x$). It is straightforward to see that if $\mu$ is
harmonic, then $\overline{\mu}$ is invariant by
$\{\sigma_{t}\}_{t\ge 0}$, and reciprocally.

\vspace{0.2cm}

If a harmonic measure can not be written as a convex sum of
different harmonic measures it is called \textit{ergodic}. The
Random Ergodic Theorem, due to Kakutani (\cite{Kakutani}, see also
\cite{Furman}) states that if $\mu$ is an ergodic harmonic measure,
then $\overline{\mu}$ is an ergodic invariant measure of
$\{\sigma_{t}\}_{t\ge 0}$. This implies in particular that for
$\mu$-almost every point $x$ in $M$, $W_x$-almost every path $\gamma
\in \Gamma_x$ is distributed with respect to $\mu$. Thus, the
ergodic properties of harmonic measures on the foliation $\mathcal
F$ can be studied via the classical ergodic theory of one
dimensional semi-groups of transformations.

\vspace{0.2cm}

Most of this work deals with 
the \textit{Lyapunov exponent} of a harmonic measure for a
transversely conformal foliation $\mathcal F$ of class $C^1$ (see
\cite{Candel,Deroin}). Let $|\cdot|$ be a transverse metric. Then if
$\gamma$ is a continuous path of $\Gamma_x$ starting at a point $x$,
consider the Lyapunov exponent
\[ \lambda(\gamma) := \lim_{t\rightarrow \infty}
\frac{1}{t} \log |Dh_{\gamma|_{[0,t]}}|,\]
when it is defined. Recall that $h_{\gamma|_{[0,t]}}$ is the
holonomy map from a transversal $T_{\gamma(0)}$ passing through
$x=\gamma(0)$ to a transversal $T_{\gamma(t)}$ passing through
$\gamma(t)$. By the Birkhoff Ergodic Theorem and the Random Ergodic
Theorem, if $\mu$ is an ergodic harmonic measure, then for
$\mu$-almost every point $x$, and $W_x$-almost every Brownian path
$\gamma$ starting at $x$, the Lyapunov exponent of $\gamma$
converges to a number depending only on $\mu$; we call it the
Lyapunov exponent of the measure $\mu$ and denote it by
$\lambda(\mu)$.

\vspace{0.2cm}

We begin by the study of the case where the Lyapunov exponent is negative.

\begin{main}\label{t:A}
Let $\mathcal F$ be a transversely conformal foliation of class
$C^1$ of a compact manifold, and suppose that on a minimal set
$\mathcal M$ is supported a harmonic ergodic measure $\mu$ with
negative Lyapunov exponent. Then the following properties are
satisfied:
\begin{itemize}
\item \textbf{Contraction.}
Let $\alpha$, $0< \alpha <|\lambda(\mu)|$ be chosen. Then for any
$x\in \mathcal M$, and almost every Brownian path $\gamma\in
\Gamma_x$, there exist a transversal $T_{\gamma}$ at $x$ and a
constant $C_{\gamma}>0$ such that for every $t>0$, the holonomy map
$h_{\gamma|_{[0,t]}}$ is defined on $T_{\gamma}$, and
\[ |h_{\gamma|_{[0,t]}}(T_{\gamma}) | \leq C_{\gamma}\exp (-\alpha t).\]
\item \textbf{Unique ergodicity.}
For any point $x\in\mcM$, almost every Brownian path starting at $x$
is distributed with respect to $\mu$. Thus $\mu$ is the unique
harmonic measure on $\mathcal M$.
\item \textbf{Diffusion.}
The diffusions $D^t f$ of a continuous function $f:\mathcal
M\rightarrow {\bbR}$ converge uniformly to the constant function
$\int_{\mathcal M} f d\mu$.
\item \textbf{Attraction.}
Suppose $\mcM\neq M$, and let $p_{\mcM}(x)$ be the probability that
a Brownian path starting at $x$ tends to $\mathcal M$, is
distributed with respect to $\mu$, and contracts exponentially a
transversal at~$x$. Then $p_{\mcM}$ is lower semi-continuous and
leafwise harmonic. In particular, $p_{\mcM}$ is bounded from below
be a positive constant in some neighborhood of~$\mcM$.
\end{itemize}
\end{main}

Theorem~\ref{t:A} is proved in section \ref{s:Negative}. The idea is
the following. A lemma of contraction (Lemma~\ref{LC}) implies,
together with the fact that the Lyapunov exponent is negative, that
for $\mu$-almost every point $x$, almost every Brownian path
starting at $x$ contracts a transverse ball exponentially. We prove
that if $\overline{x}$ is a point which is close to such a point
$x$, then there is a \textit{similarity} between the Brownian motions on the
leaf $L_{\overline{x}}$ and on $L_x$. This comes from the fact that
for a lot of Brownian pathes on $L_x$, the leaves approach each other
exponentially. All the properties annouced in the Theorem~\ref{t:A}
are deduced from this property of similarity.

\begin{remark} Theorem~\ref{t:A} is also true if the foliation is singular, but the
minimal set does not contain any singularity. In particular, our
result applies for singular holomorphic foliations on complex
compact surfaces. For instance, we prove
that if $\mathcal M$ is a minimal subset of a holomorphic
foliation of the complex projective plane, then on $\mathcal M$ is
supported a unique harmonic measure and the
Lyapunov exponent is negative. This has been recently proved by Fornaess and
Sibony for a lamination by holomorphic curves of class $C^1$ 
contained in the complex projective plane~\cite{Fornaess-Sibony}.
\end{remark}

In section \ref{s:Symmetric}, we prove that there is a dichotomy
between the case of negative Lyapunov exponent and the case where
there exists a transversely invariant measure:

\begin{main}\label{t:B}
Let $\mathcal F$ be a transversely conformal foliation of a compact
manifold. Then on a minimal subset, either there exists a
transversely invariant measure, or the harmonic measure is unique
and the Lyapunov exponent is negative.
\end{main}

In the case of a group of diffeomorphisms, the first result of this
kind was proved by Furstenberg~\cite{Furstenberg}: if $G$ is an
irreducible subgroup of projective transformations of $\bbR P ^n$,
equipped with a probability measure of finite first moment, then it
has the contraction property. For a group of diffeomorphisms of
class $C^1$ of an arbitrary compact manifold, the dichotomy ``there
exists a measure which is invariant by all the elements of the group
or there is a stationary measure with negative sum of the Lyapunov
exponents'' was proved by Baxendale~\cite{Baxendale}.


\vspace{0.2cm}

From Theorem \ref{t:B} and a Theorem of Candel~\cite{Candel}, we obtain the following:

\begin{corollary} A minimal subset of
a transversely conformal foliation of a compact manifold (if the
codimension is one, the foliation is supposed of class $C^1$)
carries an invariant measure, or there is a loop contained in a leaf
with hyperbolic holonomy. \label{hyperbolic holonomy} \end{corollary}

Our main theorem is a consequence of Theorem~\ref{t:A} and~\ref{t:B}.
At the end of section \ref{s:Symmetric}, we prove the main theorem.

\vspace{0.2cm}

In \cite{Candel}, Candel extends Garnett's theory to the case of a non symmetric Laplace
operator on the leaves of a foliation. For this kind of processes, Theorem A is still valid:
negative Lyapunov exponent implies unique ergodicity, contraction
etc. However, in the non-symmetric case Lyapunov exponent can be
positive, the dynamics not uniquely ergodic, even if the foliation
is minimal. Such an example is presented in Section~\ref{s:counter-example}.

\vspace{0.2cm}

The last part is a tentative to prove unique-ergodicity for a
Laplace operator whose drift vector field preserves the Riemannian volume.
We prove this when the foliation together with the Laplace operator are \textit{similar}.
A foliation equipped with a Laplace operator on the leaves is
similar if there exists a transverse foliation which leaves the

Laplace operator invariant.

\begin{main} \label{t:C} Let $(\mathcal F,\Delta)$ be
a similar codimension one foliation of a compact manifold, which is
transversely continuous. Suppose that the operator $\Delta$ is
obtained by drifting the Laplacian of a Riemannian metric by a
vector field that preserves the volume. Then on every minimal subset
is supported a unique harmonic measure. Moreover, if $\Delta$ is
symmetric, every ergodic harmonic measure is supported on a
minimal set.\end{main}

The idea of the proof of Theorem~\ref{t:C} is, as in
Theorem~\ref{t:A}, based on the fact that the leaves through close
points stay close in a lot of ``directions''. This idea is due to
Thurston (\cite{Thurston}, see also~\cite{Calegari,Fenley}). Because
we already know that the Brownian motion on different leaves are
similar, such a property implies unique-ergodicity. This property is
proved by constructing a harmonic transverse distance and to use the
Martingale Theorem. This is done in section~\ref{s:Similar}.

\subsection{Other examples of dynamical systems}
Following the lines of the ``Sullivan's dictionary'', our results
are still valid for discrete conformal dynamical systems. The idea
is that, instead of considering a foliation by smooth manifolds, one
can consider foliations by \textit{graphs}.

\vspace{0.2cm}

Let $M$ be a compact manifold together with a conformal structure of
class $C^1$, and $\Gamma$ be a finitely generated pseudo-group of
conformal transformations of $M$. Then for any $x$ in $M$, we denote
by $O(x)$ the orbit of $x$ under the action of the elements of
$\Gamma$. Given a symmetric system of generators of $\Gamma$, we
consider a distance on every orbit $O(x)$: the distance between $x$
and a point $y$ in $O(x)$ is the minimal number of elements of the
system of generators that is necessary to map $x$ to~$y$.


\vspace{0.2cm}

Consider a family $\{\mu_x\}_{x\in M}$ of probability measures on
$M$ whose support is the orbit of $x$. We ask that for any element
$\gamma$ of the pseudo-group $G$, the function
\[  x\in \mathrm{dom}(\gamma) \mapsto \mu_x(\gamma(x))\in (0,1)\]
is H\"older. The \textit{diffusion} operator acts on the space of
continuous functions by the formula
\[  D f(x) = \int_{O_x} f d\mu_x ,\]
for any $f\in C^0(M)$ and any $x\in M$. An invariant measure by the
diffusion semi-group always exists, and is called a
\textit{stationary} measure.

\vspace{0.2cm}

Associated to this diffusion process on the leaves, we consider the
Markov process induced by the measures $\mu_x$ on the orbits of
$\Gamma$, which is the discrete analog of the Brownian motion on the
leaves of a foliation. We define the Lyapunov exponent when the
following finiteness hypothesis holds:
\[  \forall x\in M,\quad \forall \alpha >0, \quad  \int_{O_x} \exp(\alpha d(x,y)) d\mu_x(y) <\infty.\]

\vspace{0.2cm}

Our results extends to the context of a pseudo-group with the
following analogies:
\begin{itemize}
\item the leaves together with the Brownian motion process is replaced by the orbit
of a point together with the Markovian process with transition
probabilities $p(x,y)=\sum_{\gamma(x)=y} \mu_x(\gamma(x))$.
\item the symmetry condition in the discrete case means that there exists a
H\"older function $v: M\rightarrow \bbR $ such that for any points
$x,y\in M$:
\[   p(x,y) \exp(v(x)) = p(y,x) \exp (v(y)).\]
\end{itemize}

 The important examples are the pseudo-group generated by the action of
a group on ${\bf S}^n$ by conformal transformations, or the action
of an affine conformal correspondence on a torus ${\bf T}^n$.

\newpage

\section{Negative Lyapunov exponent}\label{s:Negative}

In this section, we are going to prove Theorem \ref{t:A} for a
codimension one foliation of class $C^1$. In this case, we can find
a transversal foliation of dimension one and of class $C^1$, and
this will simplify the proof. Such a foliation does not necessarily
exist in the case of higher codimension, but at the appendix
(Section~\ref{ss:higher_codimension}) we explain how to adapt the
proof. Finally, we may suppose $\mcF$ to be transversely orientable
(and we do so from this moment), as it is always true up to a
2-folded cover, and passing a finite cover does not change our
results.

\vspace{0.2cm}

Let $\mathcal F$ be a codimension one foliation of class $C^1$ of a
compact manifold $M$, and $\mathcal G$ a transverse foliation of
class $C^1$. Consider some minimal subset $\mcM\subset M$. Let us
suppose that for some ergodic harmonic measure $\mu$ supported on
$\mcM$, we have $\lambda(\mu)<0$.

\subsection{Contraction}
The goal of this paragraph is to prove that for generic Brownian
paths, the holonomy contracts a transverse interval exponentially,
which is true infinitesimally:

\begin{proposition}\label{p:contr}
    Suppose that there exists an ergodic harmonic
    measure $\mu$ on $M$, such that $\lambda(\mu)<0$. Let
    $\alpha>0, \alpha<|\lambda(\mu)|$. Then
    for $\mu$--almost every $x\in M$ and for $W_x$-almost every
    $\gamma\in \Gamma_x$, there exists a
    transversal neighborhood $I$ and a constant $C>0$, such that
$$
\forall t>0 \quad \left|h_{\gamma|_{[0,t]}}(I)\right| < C |I|
e^{-\alpha t}.
$$
    In particular, all the holonomy maps $h_{\gamma|_{[0,t]}}$ are defined in
    the same transversal neighborhood $I$ and, as $t\to\infty$, this
    neighborhood is exponentially contracted.
\end{proposition}

To this end, we have to connect the derivatives of holonomy maps in
one point $x$ and the diameter of the images of the transversal
neighborhood $I$. The following Contraction Lemma is in the
``folklore". For the completeness of the text, we present here both
its statement and proof. This lemma extends the Distortion Lemmas,
used in $C^2$ case by Schwartz~\cite{Schwartz}, Denjoy~\cite{Denjoy}
and Sacksteder~\cite{Sacksteder}, in $C^{1+\tau}$ by
Sullivan~\cite{Su2} and Hurder~\cite{Hurder}.

\begin{lemma}[Contraction Lemma]\label{LC}
Let $x_0,x_1,\dots\in\bbR$ be points in $\bbR$, $I_j\subset \bbR,
\quad I_j=U_{\varepsilon}(x_j)$ be their $\varepsilon$-neighborhoods. Let
$$
h_j:I_j \to \bbR, \quad h_j(x_j)=x_{j+1}, \quad j=0,1,2,\dots
$$
be $C^1$-diffeomorphisms onto their image. Let
$\overline{f}_j(y)=h_j(y+x_j)-x_{j+1}$ be diffeomorphisms of
$U_{\varepsilon}(0)$, and suppose, that these diffeomorphisms are
bounded in the $C^1$ topology.

Denote
$$
F_n=h_n \circ \dots \circ h_1 : I_1 \to \bbR, \quad n\in{\bf N},
$$
and suppose that
\begin{equation}\label{eq:limsup}
\limsup_{n\to\infty} \frac{1}{n} \log F_{n}'(x_0)=\lambda<0,
\end{equation}
and let $\alpha>0, \alpha<|\lambda|$. Then there exist such an
$\varepsilon_1>0$ and a constant $C$, that for any interval
$J\subset I_0$, such that $x_0\in J, \, |J|<\varepsilon_1$ all the
compositions $F_n$ are defined on $J$ and we have a bound
$$
\forall n\in\bbN \quad |F_n(J)|\le |J| \cdot C e^{-\alpha n}.
$$
\end{lemma}

\begin{pf}
Let us choose $\beta$, $\alpha<\beta<|\lambda|$. Then, the
condition~\eqref{eq:limsup} implies that the supremum
$$
C=\sup\limits_{n\in\bbN} e^{\beta n} F_{n}'(x_0)
$$
is finite. Then, for every $n$ we have $F_n'(x_0)\le C e^{\beta n}$.

Note, that due to pre-compactness property the logatithms of the
derivatives $\log h_n'$ are uniformly continuous. Thus, there exists
$\varepsilon_0>0$, such that for every $n$ and for every $y,z\in
I_n$, $|y-z|<\varepsilon_0$ we have $h_n'(y)/h_n'(z)<
e^{\beta-\alpha}$.

Now, let $\varepsilon_1=\min(\varepsilon_0,\varepsilon)/C$, and let
$J$ be an interval of length less than $\varepsilon$,
containing~$x_0$. We are going to prove the following statement: for
every $n$, the length of $|F_n(J)|<C e^{-\alpha n} |J|$. In fact, by
the mean value theorem $|F_n(J)|= F_n'(y) |J|$ for some point $y\in
J$. Now,
\begin{equation}\label{eq:derivative}
\frac{F_n'(y)}{F_n'(x)}=\prod_{j=1}^n
\frac{h_j'(y_{j-1})}{h_j'(x_{j-1})} <
(e^{\beta-\alpha})^n=e^{(\beta-\alpha)n}.
\end{equation}
Here $y_j=F_j(y)$, and the inequality
$\frac{h_j'(y_{j-1})}{h_j'(x_{j-1})} < e^{\beta-\alpha}$ is
satisfied due to the recurrence hypothesis and the choice of
$\varepsilon_0$:
$$
|x_{j-1}-y_{j-1}|\le |F_{j-1}(J)| \le C e^{-\alpha n} |J| \le C
\varepsilon_1 < \varepsilon_0.
$$
Now, from~\eqref{eq:derivative} we have:
$$
F_n'(y)\le e^{(\beta-\alpha)n} F_n'(x) \le e^{(\beta-\alpha)n} \cdot
C e^{\beta n} = C e^{-\alpha n}.
$$
Hence,
$$
|F_n(J)| = F_n'(y) \cdot |J| < C e^{-\alpha n} |J|.
$$
This proves the recurrence step, and thus the entire lemma.
\end{pf}

\vspace{0.2cm}

Now we are going to apply the Contraction Lemma to
the Brownian motion on the leaves. In order to do that, we would
like to decompose the holonomy map in time $t$ as a composition of
some number $n$ of maps, forming a $\Diff^{1}$-pre-compact set, with
the quotient $n/t$ being bounded from above and from below.

\vspace{0.2cm}

Suppose $\delta>0$ be given. Then to any
trajectory~$\gamma\in\Gamma_x$, we associate a sequence of points
$(x_n)=(\widetilde{\gamma}({n\delta}))$ on the universal cover
$\widetilde{\mcF}_x$, where the path $\widetilde{\gamma}$ is the
covering path for the path $\gamma$, and a sequence of numbers
$k_n=[\dist(x_{n-1},x_n)]+1$ (here $[z]$ denotes the integer part
of~$z$). Let us divide a segment of shortest geodesic line, joining
$x_{n-1}$ and $x_n$, in $k_n$ equal parts; let us denote the
vertices of this partition by $y_n^0=x_{n-1}, y_n^1,\dots,
y_n^{k_n-1}, y_n^{k_n}=x_n$. Then, the holonomy map $h_{x_0 x_n}$
between $x_0$ and $x_n$ can be written as
\begin{equation*}
    h_{x_0 x_n} =h_{x_{n-1} x_n} \circ h_{x_{n-2} x_{n-1}} \dots
    \circ h_{x_0 x_1}
\end{equation*}
and thus as
\begin{equation}\label{eq:hrep}
    h_{x_0 x_n} = (h_{y_n^{k_n-1} y_n^{k_n}} \circ \dots \circ
    h_{y_n^0 y_n^1}) \circ \dots \circ
    (h_{y_1^{k_n-1} y_1^{k_n}} \circ \dots \circ h_{y_1^0 y_1^1}).
\end{equation}

The total number of the maps in the right hand side
of~\eqref{eq:hrep} is $K_n=k_1+\dots+k_n$. The following lemma is
some form of the Large Numbers Law. Though it is rather clear that
such statement should take place, its rigorous proof is rather long,
and we have put it in Section~\ref{s:TAP}.

\begin{lemma}\label{p:discrcomp}
    There exists a constant $c>0$, such that
    for any $x\in M$ for $W_x$--almost every path $\gamma\in\Gamma_x$
    we have $K_n/n < c$ for all $n$ sufficiently big.
\end{lemma}

Thus, the discretization of the Brownian motion is
``quasi-preserving'' the time: the number of terms in the right-hand
side of the representation~\eqref{eq:hrep} is comparable with the
time passed.

\vspace{0.2cm}

\begin{pfof}{Proposition \ref{p:contr}}
Suppose that a point $x\in M$ is such that for almost every path
$\gamma\in\Gamma_x$ we have
\begin{equation}\label{limLyap}
\lim\limits_{t\to\infty} \frac{1}{t} \log h_{\gamma|_{[0,t]}}'(x) =
\lambda(\mu)<0.
\end{equation}
Let us show that for this point the conclusion of
Proposition~\ref{p:contr} holds. This will prove the Proposition~---
as $\mu$--almost every point $x$ satisfies~\eqref{limLyap}.

Recall that $\alpha<|\lambda|$. As it follows from~\eqref{limLyap},
for $W_x$--almost every path $\gamma\in\Gamma_x$ there exists a
constant $C_0>0$, such that for every $t>0$
\begin{equation}\label{limalpha}
    h_{\gamma|_{[0,t]}}'(x) < C_0 e^{-\alpha t}.
\end{equation}
Let us consider a discretization $(x_n,k_n)$ of such path~$\gamma$.
As it follows from Proposition~\ref{p:discrcomp}, for almost every
path~$\gamma$ we have also
\begin{equation}\label{discr_typ}
\exists N: \forall n>N \quad \frac{K_n}{n}<c.
\end{equation}
Suppose that for $\gamma$ both~\eqref{limalpha}
and~\eqref{discr_typ} are satisfied. Choosing some $c'>c$, we may
suppose that for all $n\in\bbN$ we have $K_n<c'n$. As it follows
from~\eqref{eq:hrep}, for every $n\in\bbN$ the holonomy map $h_{x_0
x_n}$ can be written as a composition of $K_n<c'n$ maps $h_{y_l^j
y_l^{j+1}}$, each one being a holonomy between two points at the
distance at most~$1$. The set of holonomy maps along paths of length
at most 1 is pre-compact (it is a continuous image of a compact
set), and the derivative of such composition at the point $x$ is
less then $C_0 e^{-\alpha n}< C_0 e^{-\alpha K_n/c'}$, thus, we
still have an exponential decrease of derivatives at~$x$ with
respect to the number of maps. The application of the Lemma~\ref{LC}
concludes the proof.
\end{pfof}

\subsection{Similarity of the Brownian motions on different leaves.}
\label{ss:similarity}
Suppose that the point $x$ is typical in the sense of
Proposition~\ref{p:contr}, and that almost every Brownian path
starting at $x$ is distributed with respect to $\mu$. Then there
exist a transversal interval~$I$, and constants $C_0$, $\alpha>0$,
such that the set
\begin{equation*}
    E_x=\left\{\gamma\in\Gamma_x \mid \forall t \quad
    \left| h_{\gamma|_{[0,t]}}(I) \right| < C_0 \exp(-\alpha
    t) |I| \right\}
\end{equation*}
has positive Wiener measure: $W_x(E_x)>0$ (in fact, this probability
can be made arbitrary close to~$1$ by choosing sufficiently
small~$I$ and sufficiently big~$C_0$). Let us fix such $I$
and~$C_0$.

If $\bar{x}$ is a point close to $x$, consider the set $E_{\bar{x}}$
of Brownian paths $\bar{\gamma}\in \Gamma_{\bar{x}}$ which are
exponentially asymptotic to a path $\gamma$ of $\Gamma_{\bar{x}}$:
\[ \forall t\geq 0, \quad \quad d(\gamma(t), \bar{\gamma}(t) )
< C_0 \exp(-\alpha t) |I|.\] For instance, a path of $E_{\bar{x}}$
can be constructed from a path of $E_x$ by following the foliation
$\mathcal G$.

\begin{lemma}[Similarity of the Brownian motions]\label{l:positive}
There exists a neighborhood $U$ of $x$ such that for any $\bar{x}\in
U$, the Wiener measure of $E_{\bar{x}}$ is positive. Moreover, it
can be made arbitrarily close to $1$ by choosing sufficiently small
$I$, large $C_0$ and small $U$. Finally, $W_{\bar{x}}$-almost every
path $\gamma\in E_{\bar{x}}$ is distributed with respect to~$\mu$.
\end{lemma}

We are going to outline the proof of this Lemma. However, a formal
realization of the ideas faces some difficulties and becomes very
technical. Thus, we have postponed it until Appendix (section
\ref{s:TAP}).

\vspace{0.2cm}

First, we observe that the lemma is simple to prove when the
foliation is similar, meaning that the foliation $\mathcal G$
preserves the Laplace operator. In this case, the $\mcG$-along
holonomy preserves the metric and thus translates the Brownian
motion on an initial leaf to the Brownian motion of the range leaf.
In particular, the probabilities $W_x(E_x)$ and
$W_{\bar{x}}(E_{\bar{x}})$ coincide if the points $x$ and $\bar{x}$
are in the same $\mathcal G$-leaf. On the other hand, if we consider
Brownian motions starting at a point $y$ and at some point
$\bar{y}\in\mcF_y$ close enough to~$y$, the distribution of their
values at the moment $\delta$ will be absolutely continuous with
respect to each other, and the density will be close to~1. Thus, the
probabilities of any tail-type properties (in particular, of
belonging to $E_{y}$ and $E_{\bar{y}}$, being distributed with
respect to $\mu$, etc.) are close enough. The $W_x(E_x)>0$ property
implies the same for the points close enough on the same leaf of
$\mcF$ and for the points close enough on the same leaf of $\mcG$.
Thus, this property is satisfied in some neighborhood of~$x$ (as
$\mcF$ and $\mcG$ are transversal). This proves the lemma for
similar foliations.

\vspace{0.2cm}

In the general case, the $\mcF$-leafwise implication is still valid.
Unfortunately, it is much more difficult to prove the $\mcG$-along
implication, when the Riemannian metric is not invariant by the
foliation $\mathcal G$. Let us suppose that $\bar{x}\in \mcG_x$.
Denote by $\Phi_{x,\bar{x}}:E_x \rightarrow E_{\bar{x}}$ the
holonomy map along the transversal foliation~$\mcG$. Except for
similar foliations, $\Phi_{x,\bar{x}}$ does not translate the Wiener
measure $W_x$ on $E_x$ to a measure which is absolutely continuous
with respect to the Wiener measure $W_{\bar{x}}$ on $E_{\bar{x}}$.
This effect comes from small movements~--- typical path of the
Brownian motion, being considered on arbitrary small interval of
time, allows to reconstruct the Riemannian metrics (on its support).
Thus, we are going to pass from the Brownian paths to their
discretization~--- as it was already done in the previous paragraph.

\vspace{0.2cm}

Let $\delta >0$ be given. Let us define the \textit{discretization
map} $F^{\delta}: E_{x}\to (\widetilde{\mcF}_{x})^{\infty}$ as
\begin{equation}
    F(\gamma)=\{\widetilde{\gamma}(n\delta)\}_{n=0}^{\infty}.
    \label{eq:discr_def}
\end{equation}
Denote by $E_{x}^{\delta}$ the image of $E_x$ under $F^{\delta}$,
and by $W_x^{\delta}$ measure on $E_x^{\delta}$ that is the image of
the Wiener measure on $\Gamma_x$, restricted on $E_x$, under
$F^{\delta}$. We claim that if $\bar{x}\in\mcG_x$ is sufficiently
close to $x$, then the induced map
\[\Phi_{x,\bar{x}}^{\delta}: (E_x^{\delta}, W_x^{\delta})
\rightarrow (E_{\bar{x}}^{\delta}, W_{\bar{x}}^{\delta}) \] is
absolutely continuous and its Radon-Nykodym derivative is uniformly
close to $1$ on a set of a large measure. This comes from the fact
that for a finite number of steps of discretization the density can
be written explicitly: the density is the product of quotients of
heat kernels on these two leaves. The step of discretization is
constant and equals $\delta$, and the leaves approach each other
along the trajectories of $E_x$ exponentially. Thus, it is natural
to expect that such product would converge. Moreover, the closer
initial points $x$ and $\bar{x}$ are, the closer to~$1$ will be the
product. These facts imply that the probabilities $W_x(E_x)$ and
$W_{\bar{x}}(E_{\bar{x}})$ are sufficiently close and thus imply the
Lemma; their rigorous proof can be found in the appendix
(Section~\ref{s:TAP}).

To simplify the notations, we note $\sigma=\sigma_{\delta}:\Gamma\to
\Gamma$:
$$
\forall t\ge 0 \quad \sigma(\gamma)(t)=\gamma(t+\delta).
$$

\begin{corollary}[Similarity of Brownian motions]\label{c:similarity}
For any $y\in\mcM$ for $W_y$-almost every path $\gamma\in\Gamma_y$
there exists $n$, such that $\gamma(n\delta)\in U$,
$\sigma^n(\gamma)\in E_{\gamma(n\delta)}$, and $\gamma$ is
distributed with respect to~$\mu$.
\end{corollary}

\begin{pf} Let us consider the map $D:\Gamma\to\Gamma$, defined in the
following ``algorithmic'' way:
$$
D=D_1\circ D_2,
$$
where $D_2(\gamma)=\sigma^n(\gamma)$,
$n=\inf\{j\mid\gamma(j\delta)\in U\}$,
$$
D_{1}(\gamma)=\left\{
\begin{array}{ll}
\text{STOP}, & \text{if } \gamma\in E_{\gamma(0)},
\\
\sigma^n(\gamma), & n=\inf\{j\mid \,
|h_{\gamma|_{[0,j\delta]}}(I)|\ge C e^{-\alpha \cdot j\delta} |I|\}.
\end{array}
\right.
$$
Note, that $D_2(\gamma)$ is defined on $W_z$-almost every
$\gamma\in\Gamma_z$ for every $z\in\mcM$ due to the minimality of
$\mcM$ and due to the Markovian property of the Brownian motion: on
every step of discretization, the probability of hitting $U$ is
positive and bounded from below; thus the probability of never
hitting $U$ is 0.

Now, note that for every $z$ the probability of $D$ returning
``STOP'' on $\gamma\in\Gamma_z$ is bounded from below due to
Lemma~\ref{l:positive}. Also, due to the Markovian property for
every $z,z'$ the conditional distribution of $D(\gamma)$,
$\gamma\in\Gamma_z$ with respect to the condition $D(\gamma)(0)=z'$,
coincides with $W_{z'}$. Thus, the probability of $D$ not stopping
on $\gamma$ in $k$ iterations tends to zero exponentially. In
particular, $W_z$-almost surely some iteration $D^k(\gamma)$ stops,
which means, that $D_2(D^{k-1}(\gamma)\in E_{\bar{z}}$, where
$\bar{z}=D_2(D^{k-1}(\gamma))(0)$. Together with the definitions of
$D_1$ and $D_2$, this proves the first part of the corollary.

Finally, recall that for any $n$ the conditional distribution of
$\sigma^n(\gamma)$ with respect to the condition $\gamma(n\delta)=z$
coincides with~$W_z$. Lemma~\ref{l:positive} states, that for every
$z\in U$ almost every path $\gamma\in E_z$ is distributed with
respect to $\mu$, hence, the constructed path $\sigma^n(\gamma)$
almost surely is distributed with respect to~$\mu$. A finite shift
can not change the asymptotic distribution of a path, and so the
corollary is proven.
\end{pf}

\vspace{0.2cm}

Before going into the applications of the preceeding results, we
would like to make a disgression that may clarify the meaning of
``similarity of the Brownian motion on the leaves". Recall that
every Riemannian manifold $N$ of bounded geometry has a boundary
associated to the Brownian motion process on it: this is the Poisson
boundary $P(N)$. We recall the following facts, that characterize
the Poisson boundary (see \cite{Kaimanovich}). For every $x$ of $N$
there is a canonical projection $\pi _x: \Gamma_x \rightarrow P(N)$.
The family of probability measures $\nu_x = (\pi_x) _* W_x$ depends
harmonically on $N$, in the sense that for any bounded measurable
function $f$ on $P(N)$, the function
\[ P(f) (x) = \int_{P(N)} f d\nu_x\]
is a bounded harmonic function on $N$. Moreover, every bounded
harmonic function on $N$ is obtained in this way.

\vspace{0.2cm}

There is a following question, which (even if the answer is
negative), in our opinion, clarifies the proof of
Lemma~\ref{l:positive}, giving the good general idea.

\begin{question}
Given two leaves $\mathcal F_x$ and $\mathcal F_{\bar{x}}$, does the
argument of the proof of Lemma~\ref{l:positive} allow to identify
large parts of the Poisson boundaries of their universal
covers~$\widetilde{\mcF}_{x}$ and $\widetilde{\mcF}_{\bar{x}}$,
corresponding to the couple of ``directions'' in which the leaves
are converging exponentially to each other?
\end{question}

\vspace{0.2cm}

In the case of a foliation by hyperbolic surfaces of a $3$-manifold,
Thurston has constructed the ``circle at infinity'' (see \cite{Calegari,Fenley}).
A leaf of the universal cover of such a foliation is isometric to
the hyperbolic plane, and its boundary (as a Gromov hyperbolic
space) is a topological circle. Thurston has proved that there is a
natural topological identification of the circle at infinity of the
leaves of the foliation on the universal cover.

\begin{question} The boundary (as a Gromov hyperbolic space)
of the hyperbolic plane is also the Poisson boundary. Is it true
that the topological identifications of the boundaries of Thusrton's
theorem preserve the structure of the Poisson boundary as
well?\end{question}

\subsection{Proof of Theorem \ref{t:A}}

\subsubsection{Contraction property}
Let $x\in \mathcal M$. Then by Corollary \ref{c:similarity}, for
almost every Brownian path $\gamma$ starting at $x$, there exists
$n$ such that $\gamma(n\delta)$ belongs to $U$, and $\sigma^n \gamma
\in E_{\gamma(n\delta)}$. Thus starting from the time~$n\delta$, the
path $\gamma$
contracts a transverse interval exponentially, with
exponent~$\alpha$. The contraction property is proved.

\subsubsection{Unique ergodicity}
Let $x\in \mathcal M$. Recall first that almost every Brownian path
starting at $x$ is distributed with respect to $\mu$, as it is
claimed by Corollary~\ref{c:similarity}.

Let us prove that $\mu$ is the only harmonic measure supported on
$\mathcal M$. Let $\mu'$ be an ergodic harmonic measure. Then for
$\mu'$-almost every point $y$, $W_y$-almost every Brownian path
$\gamma \in W_y$ is distributed with respect to $\mu'$. Fix one of
these points $y$. By the preceeding argument, $W_y$-almost every
Brownian path $\gamma$ is also distributed with respect to~$\mu$.
Thus $\mu'= \mu$ and $\mu$ is the unique harmonic measure supported
on $\mathcal M$.

\subsubsection{Attraction}
The lower semicontinuity of the function $p_{\mcM}$ is immediately
implied by the proof of Lemma~\ref{l:positive}. Now,
$p_{\mcM}|_{\mcM}=1$ due to Corollary~\ref{c:similarity}. Thus,
$p_{\mcM}$ is bounded from below by a positive constant in some
neighborhood of~$\mcM$.

Now, due to the Markovian property of the Brownian motion, for any
initial point $x$ the process $p_{\mcM}(\gamma(t))$ is a martingale.
Due to Ito formula, $p_{\mcM}$ is leafwise harmonic.

\vspace{0.2cm}

The Theorem~\ref{t:A} is proven unless for the diffusion part, which
is absolutely analogous to the proof of the diffusion part of the
Main Theorem, see paragraph~\ref{ss:main theorem}.

\begin{remark}\label{r:typical}
In fact, for any point $y\in M$ almost all the paths
$\gamma\in\Gamma_y$, tending to $\mcM$ (if such paths exist), are
distributed with respect to~$\mcM$. In particular,
$$
p_{\mcM}(x)=W_x(\{\gamma\in\Gamma_x\mid \gamma\to\mcM\})
$$

Let us prove this. For a trajectory starting at sufficiently small
distance from~$\mcM$, with the probability close to~1 this
trajectory will be distributed with respect to~$\mu$. Decompose
trajectories starting at $y$ and tending to~$\mcM$ in two parts:
finite part before arriving close to~$\mcM$ and infinite afterwards.
The Markovian property implies, that the probability of the
trajectory being distributed with respect to~$\mu$ is arbitrary
close to total probability of tending to~$\mcM$. Thus, as the
distance of decomposition can be chosen arbitrary small, almost
every trajectory, tending to~$\mcM$, is distributed with respect
to~$\mu$.
\end{remark}

Denote by $Attr(\mcM)$ the basin of attraction of~$\mcM$: this is
the union of the leaves whose closure contains $\mathcal M$. Now, it
is rather clear that the function $p_{\mcM}$ is positive exactly in
the points of $Attr(\mcM)$.

\begin{corollary}\label{c:unique_neighborhood}
Any harmonic ergodic measure different from $\mu$ has a support
disjoint from~$Attr(\mcM)$.
\end{corollary}

\begin{pf}
A trajectory starting at $y\in Attr(\mcM)$ tends to~$\mcM$ and is
distributed with respect to~$\mu$ with positive probability. If
there exists another harmonic ergodic measure $\mu'$ with the
support not disjoint from $Attr(\mcM)$, for $\mu'$-almost every
point $z\in Attr(\mcM)$ almost all trajectories, starting at $z$,
would be distributed with respect to $\mu'$, not with respect
to~$\mu$. But there are no such points, which gives us the desired
contradiction.
\end{pf}

\subsection{Examples: holomorphic foliations on complex surfaces}
Let $\mathcal F$ be a singular holomorphic foliation of a compact
complex surface $S$, and let $g$ be a hermitian metric on $T\mathcal
F$. We note $\Laplace$ the Laplacian of $g$ along the leaves of
$\mathcal F$. With the use of harmonic measures, one can extend
certain notions that we have for compact holomorphic curves; for
instance the Euler characteristic. The following definition is due
to Candel (see~\cite{Ghys4}).

\begin{definition} Let $E \rightarrow S$ be a holomorphic line bundle over $S$, and $\mu$ a harmonic measure.
The \textit{Chern-Candel class} of $E$ against $\mu$ is
\[  c_1(E,\mu):= \frac{1}{2\pi} \int_S \mathrm{curvature}(|\cdot|)d\mu,\]
where $|\cdot|$ is a hermitian metric on $E$ of class $C^2$. Recall
that the curvature of a hermitian metric is
\[\mathrm{curvature}(|\cdot|):= -2\Laplace \log |s|,\]
where $s$ is a \textit{local} non vanishing holomorphic section of $E$. Because
the curvature of two different smooth hermitian metrics on $E$
differs by the Laplacian along $\mathcal F$ of a smooth function,
the Chern-Candel class of $E$ does not depend on the choice of the
hermitian metric. The \textit{Euler characteristic} of $T$ is the
first Chern class of the tangent bundle of $\mathcal
F$.\end{definition}

The following lemma expresses the Lyapunov exponent in algebraic
terms.

\begin{lemma} Let $\mathcal M $ be a closed minimal subset, which does not contain singularities of $\mathcal F$.
Let $g$ be a hermitian metric on $T\mathcal F$ and $\mu$ a harmonic
measure supported on $\mathcal M$. Then
\[  \lambda (\mu) = -\pi c_1 (N_{\mathcal F}, \mu).\]
\label{lyapunov exponent and chern candel class} \end{lemma}

\begin{pf} Let $|\cdot|$ be a smooth conformal transverse metric of the foliation. Let us compute the curvature
of the normal bundle $N_{\mathcal F}$. Let $(z,t)$ some local
coordinates where the foliation is defined by $dt=0$. The section
$\frac{\partial} {\partial t}$ induces a non vanishing holomorphic
section of the normal bundle. Thus, the curvature of $N_{\mathcal
F}$ is
\[  \mathrm{curvature}(|\cdot|) = -2\Laplace \log| \frac {\partial}{\partial t} |.\]
One gets
\[  c_1(N_{\mathcal F}, T_{g,\mu}) = -\frac{1}{\pi} \int \Delta \log| \frac {\partial}{\partial t} | d\mu,\]
and the lemma is proved by applying Lemma~\ref{l:lyapunov formula}. \end{pf}

\begin{theorem} Let $\mathcal F$ be a singular holomorphic foliation of a complex surface,
and $\mathcal M$ be an exceptional minimal subset of $\mathcal F$.
Suppose that the normal bundle of $\mathcal F$ has a metric of
positive curvature along $\mathcal M$. Then on $\mathcal M$ is
supported a unique harmonic measure of negative Lyapunov exponent.
\end{theorem}

\begin{pf} Let $\mu$ be an ergodic harmonic measure supported on
$\mathcal M$. The normal bundle of $\mathcal F$ has a metric of
positive curvature. Thus, by Lemma~\ref{lyapunov exponent and chern candel class},
the Lyapunov exponent of $\mu$ is negative.
The theorem is a corollary of Theorem \ref{t:A}.\end{pf}

\vspace{0.2cm}

Recently Fornaess-Sibony proved that on a compact lamination by holomorphic curves of class $C^1$ of ${\bbC}P^2$,
there is a unique harmonic measure~\cite{Fornaess-Sibony}. We recover this
result here for an exceptional minimal set of a singular holomorphic foliation of ${\bbC}P^2$.

\begin{corollary} An exceptional minimal subset of a singular holomorphic foliation of ${\bbC}P^2$
carries a unique harmonic measure which is of negative Lyapunov exponent.
\end{corollary}

\begin{pf} The normal bundle of a singular holomorphic foliation of ${\bbC}P^2$ is $O(d+2)$ where $d>0$
is the degree. This bundle has a metric of positive curvature
everywhere.\end{pf}

\newpage

\section{Symmetric case}\label{s:Symmetric}

A Laplace operator is called symmetric if it is the Laplace-Beltrami
operator of a Riemannian metric. In this section we consider
a foliation equipped with a Riemannian metric on the leaves, and the Laplace-Beltrami
operator $\Delta$ along the leaves.

\subsection{The dichotomy: proof of Theorem~\ref{t:B}}
In this paragraph we prove Theorem~\ref{t:B}: on a
minimal subset of a transversely conformal foliation is supported a
unique harmonic measure with negative Lyapunov exponent, or there is
a transversely invariant measure.

\vspace{0.2cm}

If there exists a harmonic measure for which the Lyapunov exponent is negative,
then Theorem~\ref{t:A} shows that this is the unique harmonic measure. Thus we will
prove that if the Lyapunov of \textit{any} harmonic measure is non negative,
then there exists a transversely invariant measure.

\vspace{0.2cm}

We use an integral formula which expresses the Lyapunov exponent,
and which has been founded in~\cite{Candel,Deroin}. Let $|\cdot|$ be
a transverse conformal metric. In a foliation box we consider a
transverse vector field $u$ which is invariant by the holonomy, and
define $\varphi:=|u|$. Note that $\varphi$ is well-defined up to the
multiplication by a leafwise constant function. Thus the function
$f=\Laplace \log \varphi$ is a well-defined continuous function on
$M$.

\begin{lemma} For any ergodic harmonic measure $\mu$ on $M$, $\lambda(\mu) = \int_M fd\mu$.\label{l:lyapunov formula} \end{lemma}

\begin{pf} For every $t\geq 0$, let $L_t:\Gamma \rightarrow {\bbR}$ be the functional defined by:
\[  L_t (\gamma) = \log |Dh_{\gamma|_{[0,t]}}|.\]
The family $\{ L_t\}_{t\geq 0}$ is a cocycle with respect to the shift semi-group, in the sense that for any $s,t\geq 0$:
\[ L_{t+s} = L_t + L_s \circ \sigma _t.\]
Thus the integrals $\int_{\Gamma} L_t d\overline{\mu}$ depends linearly of $t$.
Because by definition the limit
\[  \lambda (\gamma) = \lim _{t\rightarrow \infty} \frac{L_t(\gamma)} {t} \]
exists for $\overline{\mu}$-almost every path and is equal to the Lyapunov exponent of $\mu$, we have for every $t\geq 0$:
\[ \int_{\Gamma} L_td\overline{\mu}= t\lambda(\mu).\]

Now, let $x,y$ be two points in the universal cover $\widetilde L$
of a leaf $L$. Define the cocycle
\[ c(x,y) :=  \log h'_{x,y},\]
where $h_{x,y}$ is the holonomy between $x$ and $y$. Let
\[\lambda_t(x) := \int_{\widetilde{L}} p(x,y;t) c(x,y) \vol_g(y).\]
The work of Garnett (\cite[p. 288, Fact 1]{Garnett}) shows that
$\lambda_t$ is a continuous function on $M$. Moreover we have the
formula
\begin{equation}\label{eq:ll}
t \lambda(\mu) =\int_{\Gamma} L_t d\overline{\mu}= \int _M
\lambda_t(x) d\mu(x).
\end{equation}

Observe that $\varphi$ is a function defined up to multiplication by
a constant on $\widetilde{L}$ and that by definition
$\varphi(y)/\varphi(x)= h'_{x,y}$. Thus we get $c(x,y)= \log \varphi
(y) - \log \varphi(x) $ and taking the derivative at $t=0$:
\begin{equation}\label{eq:lD}
\left. \frac{d \lambda_t(x)}{dt}\right|_{t=0} = \Delta \log \varphi
(x).
\end{equation}
Differentiating~\eqref{eq:ll} at $t=0$ and
substituting~\eqref{eq:lD}, we obtain the desired result.
\end{pf}

We will also use another formula which expresses the conservation of mass of the diffusion semi-group.
Let $v$ be the volume form on $M$, induced by the leafwise volume form $\vol _g$ and the transverse metric.

\begin{lemma} $\int_M (\Delta \varphi + |\nabla \varphi|^2 ) v = 0$. \label{conservation} \end{lemma}

\begin{pf} Let $v_t$ be the transverse volume form induced by the transverse metric. Then, in a foliation box $B\times T$,
there is a volume form $\theta $ on $T$ such that
\[  v_t = \exp(q\varphi) \theta,\]
where $q$ is the codimension of $\mathcal F$. This is by definition of $\varphi$.

Consider a partition of unity: $1=\sum_i f_i$, where the support of
each function $f_i$ is contained in a foliation box $B_i\times T_i$.
By Green's formula:
\begin{multline*}
\int_M \Laplace f_i \, v = \int_{B_i\times T_i} (\Laplace f_i)
\exp(\varphi_i) \, \vol_g\wedge \theta_i =
\\
=\int_{B_i \times T_i} f_i (\Laplace \exp(\varphi_i)) \,
\vol_g\wedge \theta_i = \int_{B_i \times T_i} f_i (\Laplace \varphi
+ |\grad \varphi|^2) \exp(\varphi) \, \vol_g\wedge \theta_i =
\\
= \int_{M} f_i (\Laplace \varphi + |\grad \varphi|^2) \, v.
\end{multline*}
By summing these equalities and using the fact that $\Laplace 1=0$,
we obtain
\begin{multline*}
0= \int_M \Laplace 1 \, v = \sum_i \int_M \Laplace f_i \, v=
\\
=\sum_i \int_{M} f_i (\Laplace \varphi + |\grad \varphi|^2) \, v=
\int_{M} (\Laplace \varphi + |\grad \varphi|^2) \, v.
\end{multline*}
Lemma~\ref{conservation} is proved.
\end{pf}

\subsubsection{Proof of Theorem \ref{t:B} when $\mathcal M=M$}
By the integral definition of a harmonic measure (see~\cite[Lemma B,
p. 294]{Garnett}) and Lemma~\ref{l:lyapunov formula}, the fact that
the Lyapunov exponent of every harmonic measure is non negative is
trivially satisfied if the function $f$ is bounded from below by the
laplacian along the leaves of a smooth function $h$:
\[ f \geq  \Delta h.\]
Let us study this case as a relevant example. Let $v_t$ be the
transverse volume form induced by the transverse metric $|\cdot|$,
and $v'_t$ be defined by $v'_t= \exp(-qh)\, v_t$, where $q$ is the
codimension of $\mathcal F$. We claim that $v'_t$ is transversely
invariant, or, what is the same, the measure $\mu=\vol_g\wedge v'_t$
is totally invariant. By Lemma~\ref{conservation}, we have
\[\int_M (\Delta \varphi' + |\nabla \varphi'|^2 ) d\mu = 0.\]
Because $\varphi'$ is subharmonic, this implies that $\varphi'$
is locally constant along the leaves. Thus $\mu$ is a totally invariant measure.

\begin{remark}
A more geometric proof of this goes as follows: we consider the
leafwise gradient of the function $\log \varphi'$, which is
well-defined everywhere. Due to sub-harmonicity of $\varphi'$ it
dilates leafwise volume. By definition it also dilates the
transverse volume. Thus, the total volume must increase everywhere.
The only way is that the function $\varphi'$ must be leafwise
constant, meaning that the measure $\mu$ is totally invariant.
\end{remark}

Unfortunately, there exists a diffeomorphism of the circle, which is
minimal, and whose invariant transverse measure is
singular~\cite{Katok-Hasselblatt}. Therefore, we can not hope to
solve the functional inequality
\[  f\geq \Delta_{\mathcal F} h\]
in general. However, it is still possible to solve this inequality
``approximatively''. The following result is due to Ghys
\cite{Ghys2,Ghys3}, and the proof is based on the use of the
Hahn-Banach theorem; this idea goes back to the famous paper
\cite{cycle} by Sullivan on the foliation cycles.

\begin{lemma}\label{l:Lsolving}
Let $f:M\rightarrow \bbR$ be a continuous function such that
\[\int_M f \, d\mu\geq 0\]
for every harmonic measure $\mu$. Then there exists a sequence of
smooth functions $\psi_n : M\rightarrow \bbR$ such that uniformly,
\[ \liminf_{n\rightarrow \infty} f -\Delta_{\mathcal F} \psi_n  \geq 0.\]\label{aproximation}
\end{lemma}

\begin{pf}
In the Banach space $C^0(M)$, consider the closed subspace $E$ of
uniform limit of leafwise Laplacian of smooth functions, the cone
$\mathcal C$ of everywhere positive functions. Let $F = C^0(M) /E$
and $\overline{\mathcal C}\subset F$ be the closure of the image of
the cone $\mathcal C$ under the natural projection $C^0(M)
\rightarrow F$. The conclusion of Lemma~\ref{l:Lsolving} is
equivalent to the fact that the image of $f$ in $F$ is in
$\overline{\mathcal C}$. Suppose it is not the case. Then, by
Hahn-Banach separation theorem, there exists a linear functional
which is non-negative on $\overline{\mathcal C}$ and negative on the
image of $f$. Such a linear functional is by definition a harmonic
measure (due to the integral definition of harmonicity), if it is
conveniently normalized, and the corresponding Lyapunov exponent is
negative. It gives us the desired a contradiction. The lemma is
proved.
\end{pf}

To prove Theorem~\ref{t:A} in the case where $\mathcal M=M$, we consider the family of
volume forms
\[  \mu_n = \exp (-q \psi_n) \vol_{g} \wedge v_t,\]
where the functions $\psi_n$ are given by Lemma~\ref{aproximation}.
After normalizing them to probability measures and taking a
subsequence, they converge to a probability measure $\mu$ on $M$.

\begin{lemma}
The measure $\mu$ is transversely invariant.\label{l:invariance}
\end{lemma}

\begin{pf} It is more convenient to consider the family of transverse
forms $v_{t,n} = \exp(-q \psi_n) v_t$ as currents, i.e. as operators
on the space of $p$-forms along the leaves. Define:
\[C_n(\omega ) = \int_M \omega \wedge v_{t,n},\]
for every $p$-form $\omega$, where $p=\dim (\mathcal F)$. The family
of currents $\{C_n\}$ is bounded, thus after taking a subsequence
they converge to a current $C$. By construction, by choosing well
the subsequence of currents converging to $C$, we have
\[ \mu = \vol_g \wedge C.\]
We prove that $C$ is a \textit{closed} current, or equivalently that
$C_n$ is asymptotically closed as $n$ goes to infinity. Thus by
Sullivan's Theorem (\cite[Theorem I.12, p. 235]{cycle}) the measure
$\mu$ is totally invariant.

Consider a $(p-1)$-form $\alpha$. Using a partition of unity, we can
write $\alpha$ as a finite sum of $(p-1)$-forms whose support is
contained in a foliation box. Thus we will suppose that the support
of $\alpha$ is contained in a foliation box $B\times T$. Let $\theta$ be
a volume form on $T$; write
\[   v_{t,n} = \exp (\varphi_n ) \theta,\]
then
\begin{equation}\label{eq:hdensity}
\mu_n=v_{t,n}\wedge \vol_g = \exp (\varphi_n ) \theta\wedge \vol_g.
\end{equation}
We have
\[   \int_M d\alpha \wedge v_{t,n} = \int_M d\alpha \wedge \exp(\varphi_n)\theta =\]
\[ \int_M \alpha \wedge d(\exp(\varphi_n) \theta) =
\int_M (\alpha \wedge d\varphi_n)\wedge v_{t,n}.\] Thus, by Schwarz
inequality,
\begin{equation}\label{eq:Schwarz}
\left|\int_M d\alpha \wedge v_{t,n} \right| \leq c|\alpha |_{\infty}
\int_M |\nabla  \varphi_n| \mu_n \leq c|\alpha |_{\infty} \left(
\int_M |\nabla \varphi_n|^2 \mu_n \right)^{1/2},
\end{equation}
where $c$ is a constant, and $\nabla$ denotes the leafwise gradient.
Observe that $\varphi_n$ is well-defined up to addition of a
leafwise constant function, so that $\nabla \varphi_n$ is
well-defined.

By Lemma~\ref{conservation}:
\[  \int_M |\nabla \varphi_n| ^2 \mu_n = - \int_M \Delta \varphi_n \mu_n,\]
and by Lemma \ref{aproximation}, the Laplacians $\Delta \varphi_n$
verify uniformly
\[  \liminf_{n\rightarrow \infty} \Delta \varphi_n \geq 0.\]
Thus the integrals
\[  \int_M |\nabla  \varphi_n |^2 \mu_n \]
tend to $0$ when $n$ goes to infinity. By formula~\ref{eq:Schwarz},
we conclude that $C$ is closed. Hence Lemma~\ref{l:invariance} is
proved.
\end{pf}

The proof of this lemma concludes the proof of Theorem~\ref{t:B} when $\mathcal M=M$.

\begin{remark}
In general, for a codimension $q$ foliation of class $C^1$,
Ose\-le\-dets' Theorem states the existence of $q$ Lyapunov
exponents $\lambda_1(\mu)$,\dots, $\lambda_q(\mu)$ associated to any
harmonic ergodic measure~$\mu$. In fact, the method we have used in
this section proves that for any symmetric Laplace foliation of
class $C^1$: either there exists a harmonic measure $\mu$ for which
$\lambda_1(\mu)+\ldots + \lambda_q(\mu)$ is negative, or there
exists a transversely invariant measure. This result is analogue to
the one of Baxendale~\cite{Baxendale}.

However, let us mention some differences between the context of
groups and the one of foliations. The theorem of Baxendale does not
require the dynamics to be symmetric. It is interesting to note that
Baxendale's Theorem does not work for non symmetric Laplace
foliations: there exists an example of a minimal foliation (drifted
geodesic flow, see paragraph~\ref{s:counter-example}), for which
every harmonic measure has positive Lyapunov exponent.  The fact
that the probability does not depend on the point in the Theorem of
Baxendale should be interpreted in the foliation context by the
concept of similarity. There is also an example of a similar Laplace
lamination which does not verify Baxendale's Theorem (in this
setting the Lyapunov exponent should be defined by using the
transverse $2$-adic structure).
\end{remark}

\subsubsection{Proof of Theorem~\ref{t:B} in the exceptional minimal set case}

First, let us notice that $\mcF$ can not support a harmonic measure
with positive Lyapunov exponent. A general idea implying this is the
following. Consider first the case of a similar foliation. Then, a
harmonic measure on $\mcM$ induces a harmonic transverse measure
(see Definition~\ref{transverse harmonic measure}), which is
harmonic. Due to the Ito Formula, the transversal measure of the
image of a small transverse ball under the holonomy map
$h_{\gamma|_{[0,t]}}$ associated to the Brownian path
$\gamma|_{[0,t]}$ is a martingale.

Hence, its expectation at any Markovian moment equals the measure of
the initial transverse ball. On the other hand, a small transverse
ball around a typical point (in the sense of convergence of the
Lyapunov exponent) is exponentially expanded by a typical Brownian
path.

Thus, the martingale takes (at some big moments of time) large
values with a large probability. Thus, the expectation is large.
This contradicts the fact that the initial transverse ball can be
chosen arbitrarily small. Rigorous proof of this statement is
presented in Section~\ref{s:TAP}, Lemma~\ref{l:non-positive}. Thus,
all the Lyapunov exponents are equal to~$0$.

The fact that the Laplace operator is symmetric was used when we
stated that conditional measures are harmonic: if the Laplace
operator is non-symmetric, these measures are harmonic in the sense
of the adjoint operator~$\Laplace^*$, not $\Laplace$.

Now, let us continue the proof using the fact that all the Lyapunov
exponents vanish. Choose some $\varepsilon>0$. Then, the arguments
of Lemma~\ref{l:Lsolving} imply that there exists a function
$\psi_{\varepsilon}^0$, such that $-\varepsilon< f-\Laplace
\psi_{\varepsilon}^0 < \varepsilon$. By continuity, the same
inequality holds in some neighborhood $U^{\varepsilon}$ of~$\mcM$,
that we suppose to be contained in the $\varepsilon$-neighborhood
of~$\mcM$.

\begin{lemma}\label{l:psi} If $\varepsilon>0$ is small enough,
there exists a function $\psi_{\varepsilon}$, such that $\Laplace
\psi_{\varepsilon} \ge -\varepsilon$ and such that at least
($1-\varepsilon$)-part of the measure $\mu_{\varepsilon} =
e^{-\psi_{\varepsilon}^0+\psi_{\varepsilon}} \vol_g \wedge v_t$ is
concentrated in~$U^{\varepsilon}$.\end{lemma}

Once such functions are constructed for any $\varepsilon$, the proof
will be finished in the same way as the proof in a minimal case.
Namely, suppose that such functions are constructed. Let us find any
weak limit $\mu$ of a subsequence of a family
$\frac{1}{\mu_{\varepsilon}(M)} \mu_{\varepsilon}$ as
$\varepsilon\to 0$. Note that $\mu$ is supported on~$\mcM$: this
comes from the fact that $\mu_{\varepsilon}(U^{\varepsilon})\ge
1-\varepsilon$ and that $\bigcap_{\varepsilon}
U^{\varepsilon}=\mcM$. Also, note that $\mu$ is a weak limit for the
same subsequence of the (non-normalized) measures
$\mu_{\varepsilon}|_{U^{\varepsilon}}$. But these restricted
measures can be written (locally) as $e^{\varphi'} \vol_g\wedge
\theta$, where $\theta$ is a transverse measure, and $\Laplace
\varphi' = q\Laplace\varphi-q\Laplace \psi_{\varepsilon}^{0}+
\Laplace\psi_{\varepsilon} \ge  -2\varepsilon$. Passing to the limit
and making estimates as in Lemma~\ref{l:invariance}, we obtain that
$\mu$ is totally invariant. Thus the proof of the Theorem will be
complete after the proof of Lemma \ref{l:psi}.

\vspace{0.2cm}

\begin{pfof}{Lemma \ref{l:psi}}
Recall that $\varepsilon>0, \psi^{0}_{\varepsilon}$ and
$U^{\varepsilon}$ are already chosen. We choose another neighborhood
$V$ of $\mcM$, $V\subset U^{\varepsilon}$. We are going to find a
function $\psi = \psi_{\varepsilon}$ verifying Lemma \ref{l:psi}.

First, let us suppose that in $U^{\varepsilon}$ there is no other
minimal set, and thus that every leaf passing through a point in
$U^{\varepsilon}\setminus \mcM$ intersects $\partial
U^{\varepsilon}$. We are going to look for the function $\psi$ as a
solution of the Poisson equation
$$
\Laplace \psi (x)  = - \varepsilon, \quad x\in
U^{\varepsilon}\setminus V; \quad
\psi|_{\partial(U^{\varepsilon}\setminus V)}=0,
$$
extended by $0$ to the complementary of $U^{\varepsilon}\setminus
V$. The solution of this problem always exists and can be found in
the following way:
$$
\psi(x) = \varepsilon \Expect T(\gamma),
$$
where $T(\gamma)$ denotes the first intersection moment of a
Brownian path $\gamma$ with the boundary of
$U^{\varepsilon}\setminus V)$:
$$
T(\gamma)= \min\{t:\gamma(t)\in \partial (U^{\varepsilon}\setminus
V)\},
$$
and $\Expect$ is the expectation of a function on the probability
space $(\Gamma_x,W_x)$. Note that $\Laplace \psi$ equals
$-\varepsilon$ in $U^{\varepsilon}\setminus V$, and $0$ in $V$ and
in the complementary of $U^{\varepsilon}$. Moreover, on $\partial
(U^{\varepsilon}\setminus V)$, $\Laplace \psi$ is a positive
distribution. Thus, one has $\Laplace \psi \geq -\varepsilon$.

Now, let us show that for an appropriate choice of $V$ the major
part of $\mu_{\varepsilon}$ is concentrated in $U^{\varepsilon}$,
with the precise estimates of the Lemma. To do this, it suffices to
check that as $V$ tends to $\mcM$, the part of the measure
$\mu_{\varepsilon}$, concentrated in $U^{\varepsilon}$, tends to~1.
Note that the (non-normalized!) measure $\mu_\varepsilon$ of
$M\setminus U^{\varepsilon}$ does not change, so that we have to
prove that the measure of $U^{\varepsilon}$ tends to infinity. By
the monotone convergence theorem, it is equivalent to the fact that
if we let $\bar{\psi}=\lim_{V\to \mcM} \psi$ (maybe, $\bar{\psi}$
equals infinity at some points), the function $e^{\bar{\psi}}$ will
be non-integrable in $U^{\varepsilon}$. Note that the function
$\bar{\psi}$ can be written as:
$$
\bar{\psi}(x) = \varepsilon \Expect T_0(\gamma),
$$
where
$$
T_0(\gamma)= \min\{t:\gamma(t)\in \partial U^{\varepsilon}\}
$$
(if such an intersection does not occur, we define
$T_0(\gamma)=\infty$). Thus, we have to estimate the mean $\Expect
T_{0}(\gamma)$.

Note that in a neighborhood $U^{\varepsilon}$ we have $f-\Laplace
h_{\varepsilon}<\varepsilon$. Thus, for a distance $\tilde{d}$
induced by a transversal metric $e^{-h_{\varepsilon}} |\cdot|$, we
have $\Laplace\log \tilde{d}(\cdot,\mcM)<\varepsilon$ in
$U^{\varepsilon}$.

Now, let us consider a random process
$$
\xi_0(t,\gamma)=\log \tilde{d}(\gamma(t),\mcM) -\varepsilon t,
$$
and let us stop it at the moment $T_0(\gamma)$:
$$
\xi(t,\gamma)=\xi_0(\min(t,T_0(\gamma)),\gamma).
$$
Then, the Ito formula implies that $\xi(t,\gamma)$ is a
supermartingale:
\begin{multline*}
\left. \frac{\partial}{\partial s}\right|_{s=t+0}
\Expect(\xi(s,\gamma)\big|\gamma|_{[0,t]}) =
\\
=\begin{cases} (\Laplace
\log\tilde{d}(\cdot,\mcM))(\gamma(t))-\varepsilon , &
\gamma|_{[0,t]}\subset U^{\varepsilon}
\\
0, & T_0(\gamma)\ge t.
\end{cases}
\le 0
\end{multline*}

Note also, that the function $\log\tilde{d}(\cdot,\mcM)$ is
Lipschitz on the leaves, thus, the conditional second moments
$$
\Expect (r_n^2(\gamma)\big|\gamma|_{[0,t]})
$$
of the increasements
$$
r_n(\gamma)=\xi(n+1,\gamma)-\xi(t,\gamma)
$$
are bounded uniformly on~$n$ and $\gamma|_{[0,n]}$. Thus, due to the
theory of martingales, for every Markovian moment~$\tau$ with finite
expectation the expectation of $\xi$ at this moment does not exceed
its initial value.

Let us now use this process to estimate from below the
expectation~$\Expect T_0(\gamma)$. Either this expectation is
infinite (in which any lower bound is satisfied automatically). Or
it is finite, and in this case the expectation of a value of a
supermartinagle $\xi$ in a Markovian moment $T_0(\gamma)$ does not
exceed its initial value, that is
$$
\Expect \left[\log \tilde{d}(\gamma(T_0(\gamma)),\mcM) -\varepsilon
T_0(\gamma)\right] \le \log \tilde{d}(x,\mcM).
$$
The expectation in the left side can be rewritten as
$$
-\varepsilon \Expect T_0(\gamma) + \Expect\log
\tilde{d}(\gamma(T_0(\gamma)),\mcM)= - \varepsilon \Expect
T_0(\gamma) + O(1),
$$
for at the moment of exiting $U^{\varepsilon}$ the distance to
$\mcM$ is separated from~$0$. So, we have
$$
\bar{\psi} = \Expect T_0(\gamma)\ge -\frac{1}{\varepsilon} \log
\tilde{d}(x,\mcM) + C_0.
$$
This implies that
$$
e^{\bar{\psi}(x)}\ge \frac{C}{(\tilde{d}(x,\mcM))^{1/\varepsilon}}.
$$
Thus, if $\varepsilon$ is less than $1/(\codim \mcF)$, the function
$e^{\bar{\psi}}$ is non-integrable and the effect of concentration
takes place. This completes the proof under the hypothesis that
$\mcM$ is the unique minimal subset of $U^{\varepsilon}$.

To conclude the proof, we remark that if $U^{\varepsilon}$ contains
another minimal set, we can replace $\mcM$ by the closure of the
union of all the leaves, entirely contained in~$U^{\varepsilon}$,
and repeat the previous arguments.\end{pfof}

The proof of Theorem~\ref{t:B} is completed in all the cases.

\subsection{Proof of the Main Theorem}\label{ss:main theorem}
We suppose that the foliation $\mathcal F$
is transversely conformal and does not have a transversely invariant measure.
By Theorem~\ref{t:B}, on any minimal set is supported a unique harmonic measure with negative Lyapunov
exponent. Because of the attraction property, any minimal set has a neighborhood which does not contain
any other minimal set. Thus, there is
a finite number of minimal sets $\mcM_1,\dots,\mcM_k$.
Denote by $\mu_1,\dots,\mu_k$ their unique harmonic measure, and $\lambda_1,\ldots,\lambda_k$
the corresponding Lyapunov exponents. Note
that every point $x\in M$ belongs to the basin of attraction of at
least one of these sets; the reason is that the set
$$
M\setminus (\bigcup_{j=1}^k Attr(\mcM_j))
$$
is closed, consists only of entire leaves and does not contain any
minimal subset.

Let $\alpha >0$ be a real number such that $\alpha <|\lambda_j|$ for
every $j$. For every point $x\in M$ we consider the probability
\[  p_j(x)=W_x(\{\gamma\in\Gamma_x \mid \gamma(t)
\xrightarrow[t\to\infty]{} \mcM_j\}),\] that a Brownian path
starting at $x$ tends to $\mathcal M_j$. Note that almost every
Brownian path tending to $\mcM_j$ (if such path exists) is
distributed with respect to $\mu_j$, and contracts a transverse ball
at $x$ exponentially with exponent $-\alpha$ (see
Remark~\ref{r:typical}). We claim that the sum of these
probabilities is equal to~$1$; in other words, $W_x$-almost every
trajectory tends to one of the minimal sets, with the distribution
and transverse contraction properties. We show this in the following
way: for arbitrary small neighborhoods $U_1,\dots, U_k$ of
$\mcM_1,\dots,\mcM_k$ respectively, the complementary $R=M\setminus
\cup_j U_j$ is a closed set without any minimal subset, thus
containing no entire leaf. Hence, for any point $x$ of $R$ there
exists a leafwise path leading to one of the neighborhoods $U_j$;
moreover, by compactness of $M$, the length of such a path is
bounded uniformly on~$R$. Thus, for a point $x\in R$, the
probability that it lies in one of the $U_j$ at time $1$ is bounded
from below by a positive uniform constant. Hence, for any point
$x\in M$, almost every trajectory starting in $x$ meets one of the
neighborhoods~$U_j$.

To complete the proof, let us show that $\sum_j
p_j(x)>1-\varepsilon$ for any $\varepsilon>0$. To do this, let us
choose $U_j$ so close to $\mcM_j$ that for every point in $U_j$ the
probability of attracting to $\mcM_j$ with the distribution and the
transverse contraction properties
is at least $1-\varepsilon$;
it is possible due to Lemma~\ref{l:positive}. Now, let us use the
Markovian property: for any $x\in M$, almost every trajectory
$\gamma\in\Gamma_x$ meets one of the $U_j$, and for a starting point
in $U_j$ the probability of attracting to the corresponding $\mcM_j$
is at least $1-\varepsilon$. Thus, the probability of
attracting to one of the $\mcM_j$ is at least $1-\varepsilon$. As
$\varepsilon>0$ was chosen arbitrary, we have proven that almost
every trajectory tends to one of the $U_j$.

Recall that the functions $p_j$
are leafwise harmonic and lower
semicontinuous. Because their sum is equal to the constant function
$1$ the functions $p_j$ are continuous.

Thus, we have proved the Contraction, Distribution and Attraction
parts of the main theorem.
We end the proof of Theorem \ref{t:main theorem} by proving the
statement about the asymptotic behaviour of the diffusion. First, we
shall prove a weaker form. Namely, we prove that the
\textit{time-averages} of diffusions tends to the same limit:
$$
\frac{1}{T} \int_0^T D^t f \, dt \stackrel{x\in
M}{\rightrightarrows} \psi(x),
$$
where $\psi(x) = \sum_j  p_j \int fd\mu_j$.
In the case where the foliation is minimal, it is implied by unique ergodicity
(see analogous arguments in~\cite{Furstenberg-minimality}). Namely,
the value of the time-average of the diffusions at a point $x\in M$ can
be rewritten as an integral:
\begin{multline*}
\frac{1}{T} \int_0^T (D^t f)(x) \, dt = \frac{1}{T} \int_0^T \int_M
D^t f \, d\delta_x \, dt =
\\
= \frac{1}{T} \int_0^T \int_M f \, d(D^t_* \delta_x) \, dt = \int_M
f \, dm_{x,T},
\end{multline*}
where
$$
m_{x,T}=\frac{1}{T} \int_0^T (D^t_* \delta_x) \, dt
$$
is the time-average of the diffusions of the measure~$\delta_x$.
Note that due to a classical argument in ergodic theory, a weak
limit of a sequence $m_{x_n,t_n}$ with $t_n\to\infty$ is harmonic.
As there exists a unique harmonic measure $\mu$, the time averages
$m_{x,t}$ converge to $\mu$ uniformly in $x$ as $t$ tends to
infinity. Thus, the integrals of $f$ with respect to these measures
also converge uniformly to $\int_M f\, d\mu$, which implies the
desired statement.

In the case of an exceptional minimal set we notice that the
time-averages of the diffusions can be rewritten as
\begin{equation}\label{eq:time-averages}
\frac{1}{T} \int_0^T (D^t f)(x) \, dt = \int_{\Gamma_x} \left(
\frac{1}{T} \int_0^T f(\gamma(t)) \, dt \right) \, dW_x(\gamma).
\end{equation}
We know that $W_x$-almost all trajectories tend to one of the
$\mcM_j$'s and are distributed with respect to the corresponding
harmonic measure $\mu_j$. The probability that a point $x$ tends to
$\mcM_j$ is equal to $p_j(x)$. Hence, the right hand side
of~\eqref{eq:time-averages} is equal to
$$
\sum_{j=1}^k p_j(x) \int_{\mcM_j} f \, d\mu_j.
$$
Moreover, a uniform argument on $(1-\varepsilon)$-measure of
trajectories (similarity of Brownian motions) implies that this
convergence is uniform in~$x\in\mcM$.

Thus, in every case we have shown that the time-average of the
diffusions converge to the right-hand side, which we denote~$\psi$.

Now, let us finish the proof using the arguments analogous to these
of Kaimanovich~\cite{Kaimanovich}. Namely, notice that due to the
diffusion of the Brownian motion, there exists
$\varepsilon,\varepsilon'>0$ such that for any point $x$, and any
point $y\in U_{\varepsilon}(x)$, the densities $p(x,y,1)$ and
$p(x,y,2)$ are bounded from below by $\varepsilon'$. Thus one has
for every $x$ and every bounded function~$f$:
\[ |D^1f(x)-D^2f(x)| \leq 2(1-\varepsilon'\mathrm{vol}(U_{\varepsilon})(x))|f|_{\infty},\]
where $|\cdot|_{\infty}$ is the uniform norm. Thus, because the
leaves are of bounded geometry, we have
\[  || D^1 -D^2 || _{\infty} <2, \]
where $||\cdot||_{\infty}$ is the norm of operators acting on
$L^{\infty}$. The ``zero-two law''~\cite{Lin} implies that
$\|D^n-D^{n+1}\|\to 0$ as $n\to\infty$. In particular, the
time-averages of the diffusions converge if and only if the
diffusions converge themselves to the same
limit~\cite{Kaimanovich2}. Hence, the diffusions converge to the
limit we have described.

\vspace{0.2cm}

\begin{pfof}{Corollary~\ref{hyperbolic holonomy}}
Let $\mathcal F$ be a transversely conformal foliation of class
$C^1$ of a compact manifold, and $\mathcal M$ a minimal set of
$\mathcal F$. By Theorem~\ref{t:B}, either $\mathcal M$ supports a
transversely invariant measure, or a harmonic measure of negative
Lyapunov exponent. In this case, Candel has proved that there exists
a loop contained in a leaf of $\mathcal M$ with hyperbolic holonomy
(see~\cite[Theorem 8.18]{Candel}). \end{pfof}

\subsection{Examples: codimension one foliations of class $C^2$}
In the case of codimension one foliations of class $C^2$ without compact leaf, the following result
completes the Main Theorem:

\begin{proposition} Let $\mathcal F$ be a codimension one foliation of class $C^2$
of a compact manifold, without compact leaves. Then if $\mathcal F$ has a totally invariant measure $\mu$,
this measure is the unique harmonic measure. Then, for every point $x$, almost every Brownian path starting at
$x$ is distributed with respect to $\mu$, and the diffusions of a continuous function $f:M\rightarrow {\bbR}$
tend uniformly to the constant function $\int f d\mu$.
\end{proposition}

\begin{pf} By Sacksteder Theorem, the foliation $\mathcal F$ is
minimal. We first prove that the measure $\mu$ is the unique totally
invariant measure. By minimality of $\mathcal F$ and Haefliger's
argument \cite{Haefliger}, there exists a transverse circle $C$
cutting every leaf. The transversely invariant measure
(corresponding to~$\mu$) induces a measure $\theta$ on $C$,
invariant by all the holonomy maps. This measure gives us a map
$h:C\to\bbR/l\bbZ=C'$, where $l=\theta(C)$. This map semi-conjugates
the pseudo-group induced by $\mathcal F$ on $C$ to a finitely
generated group of rotations of $C'$, which we note $G$. Because
$\mathcal F$ is minimal, and the holonomy pseudo-group is finitely
generated, at least one of the rotations of $G$ is irrational. Then,
the Lebesgue measure is the unique probability measure invariant by
$G$, and thus $\mu$ is the unique totally invariant measure on
$\mathcal F$ up to multiplication by a constant.

To conclude the case when there exists a totally invariant measure
$\mu$, it suffices to show that every harmonic measure is in fact a
totally invariant measure. Observe that the group $G$ is a group of
rotations, so that the orbits of its action on $C'$ have a
polynomial growth. Hence, the same is true for the action of the
holonomy group on $C$. Thus, every leaf grows polynomially.
Kaimanovich (\cite[Corollary of Theorem~4]{Kaimanovich}) proved that
if for a harmonic measure, almost every leaf (with respect to this
measure) has subexponential growth, then this measure is totally
invariant. In our case, all the leaves have polynomial growth, hence
every harmonic measure is in fact totally invariant.

The distribution and diffusion property is implied by the fact that the
harmonic measure is unique.\end{pf}

We end the paragraph by constructing a foliation by surfaces of a
$3$-dimensional compact manifold, with two exceptional minimal sets.
The example is constructed in the following way. Let
$\mathcal F$ be an oriented codimension one foliation by oriented
surfaces with an exceptional minimal set $\mcM$
and suppose that there exists
a transverse loop $c$ which does not cut $\mcM$.
Then a neighborhood of $c$ in $M$
is diffeomorphic to a solid torus $D^2\times {\mathbb S}^1$,
the foliation $\mathcal F$ being the horizontal fibration by two
dimensional balls $D^2$. Now consider two copies $N_1$ and $N_2$ of
the exterior of $D^2 \times {\mathbb S}^1$ in $M$. These two manifolds
are foliated, and have a boundary component $\partial B \times
{\mathbb S}^1$ transverse to the foliation $\mathcal F$. The
foliation $\mathcal F$ induces the horizontal foliation by circles on
it. Observe also that $N_1$ and $N_2$ have an exceptional minimal
set in their interior. Thus, by gluing $N_1$ and $N_2$ along their
boundary by a diffeomorphism which preserves the foliation and
reverses the orientation, we construct a foliation by surfaces of a
closed manifold with two exceptional minimal sets.

Now we are given an example of such a situation. We consider a
surface $\Sigma$ of bounded topology and constant negative
curvature, with a cusp of infinite volume.
The cusp determines an interval $I$ in the boundary of the universal cover of $\Sigma$ which
has two remarcable properties. The first is
that it is invariant by the action of the geodesic $\gamma$ on
$\widetilde{\Sigma}$. The second is that it is a component of the
exterior of the limit set of $\pi_1(\Sigma)$. Now consider a compact
surface $S$ of sufficiently large genus so that there exists a
surjective morphism $\rho:\pi_1(S) \rightarrow \pi_1(\Sigma)$. We
get an action of the fundamental group of $S$ on the boundary of
$\widetilde{\Sigma}$, which leaves the limit set of $\pi_1(\Sigma)$
invariant, and for which there exists an element leaving $I$
invariant, and which acts as a translation on it.
Let $(M,\mathcal F)$
be the \textit{supension} of $\rho$:
this is the foliation induced by
a flat circle bundle over $S$ whose holonomy is smoothly conjugated
to the representation $\rho$ (see~\ref{suspension}). Then the saturated subset of the limit set of $\pi_1(\Sigma)$
is an exceptional minimal subset $\mathcal M$ of $\mathcal F$.
Now, by construction, there is a leaf $L$ which intersects twice the component $I$ of the exterior of the limit set
of $\pi_1(\Sigma)$ in $\partial{\widetilde{\Sigma}}$. By the standard Haefliger's argument, we construct
a transverse circle which does not cut $\mathcal M$. Thus applying the preceeding arguments we
construct a foliation with two minimal sets.

\subsection{A counter-example in the non symmetric case}\label{s:counter-example}
In \cite{Candel}, Candel extends Garnett's theory to the case of non
symmetric Laplace operators on a foliation. In the case of non
symmetric Laplace operators, the dichotomy ``The Lyapunov exponent
is positive or there exists a transversely invariant measure'' does
not hold anymore. In this paragraph we describe a nice
counter-example in the non symmetric case (see also \cite{Hamenstadt2}).

\vspace{0.2cm}

Consider a compact Riemannian manifold $(M,g)$ of dimension $3$ on
which there is an orthonormal frame $(H^s,V,H^u)$, for which the
vector fields $H^s,V,H^u$ verify the relations:
\[  [V,H^s]= -H^s,\quad [V,H^u] = H^u, \quad [H^u,H^s] = V.\]
Such manifolds are quotient of the universal cover of $SL(2,{\bf
R})$ by a cocompact lattice $\Gamma$. If $\Sigma$ is the quotient of
the upper-half plane ${\bf H}$ by $\Gamma$, then $M$ is naturally
identified with the unitary tangent bundle of $\Sigma$, and under
this identification $V$ is the geodesic flow of the hyperbolic
surface, $H^s$ and $H^u$ the horocycle flows. The vector fields $V$
and $H^s$ generate a foliation $\mathcal F^s$ which is the
\textit{stable} foliation of the flow $V$.

\vspace{0.2cm}

Let $g^s$ be the restriction of the metric $g$ on $\mathcal F^s$.
For any $\kappa\in {\bbR}$, consider the Laplace operator
$\Delta_{\kappa}$ defined by
\[   \Delta_{\kappa} = \Delta_{g^s} + \kappa V,\]
where $\Delta_{g^s}$ is the Laplacian of $g^s$ along the leaves of
$\mathcal F^s$. Garnett proved that for the symmetric case $\kappa
=0$, the Liouville measure $\mathrm{vol}_g$ is the unique harmonic
measure (see~\cite[Proposition 5, p.~305]{Garnett}). Note that the
Liouville measure on $M$ is also invariant by $V$, so that it is a
harmonic measure for all the Laplace operators $\Delta_{\kappa}$.

\begin{theorem} For any $\kappa$, the Lyapunov exponent of any harmonic measure $\mu$ of $(\mathcal F^s,\Delta_{\kappa})$
is $\lambda (\mu) = \kappa -1$. When $\kappa <1$ the Liouville
measure is the unique harmonic measure. When $\kappa >1$, there
exists a harmonic measure supported on every cylinder leaf (thus the
foliation is not uniquely ergodic).
\end{theorem}

\begin{pf} First, we compute the Lyapunov exponent of a harmonic measure $\mu$ of $(\mathcal F^s,\Delta_{\kappa})$.
To this end we use the formula of Lemma \ref{l:lyapunov formula}:
\[  \lambda(\mu) = \int _{M} \Delta_{\kappa} \log \varphi d\mu.\]
Consider the metric $|\cdot|$ on the normal bundle of $\mathcal F$
which is induced by $g$. We are going to compute the function
$\varphi$, which is defined up to multiplication by a constant. This
function verifies the relations
\[   V \varphi = \varphi,\quad H^s \varphi =0.\]
In the leaves we have local coordinates $z= x+iy$ with values in the
upper half-plane ${\bfH}$, such that
\[   V = y\frac{\partial}{\partial y},\quad H^s = y\frac{\partial} {\partial x}.\]
(These coordinates are well defined up to an affine transformation
of the upper half-plane). In these coordinates, $\varphi = y$ up to
a multiplicative constant. The metric $g^s$ and the Laplacian
$\Delta_{g^s}$ are expressed by
\[ g^s  = \frac{dx^2 + dy^2}{y^2},\quad \Delta_{g^s}= y^2 ( \frac{\partial^2}{\partial x ^2}+\frac{\partial^2}{\partial y^2}).\]
Thus, we have $\Delta_{\kappa} \varphi = \kappa - 1$ identically,
and the formula $\lambda (\mu) = \kappa -1$ follows. In particular,
when $\kappa<1$ the only harmonic measure is the Liouville measure,
because of Theorem \ref{t:A}.

Now, let us suppose that $\kappa>1$. We shall prove that every
cylinder leaf supports a harmonic measure. Observe that these leaves
are those containing a periodic orbit of the vector field $V$. Let
$L$ be such leaf, $\gamma_0$ be the closed orbit of $V$ in~$L$.
Then, its universal cover is the hyperbolic plane, for which we
choose the upper half-plane model~${\bfH}$. Observe that we have the
canonical coordinates up to an affine transformation, constructed
before. Without loss of generality, we may suppose that the geodesic
in ${\bfH}$, corresponding to $\gamma_0$, is the vertical geodesic
$x=0$ going upwards. Denote by $A$ the length of $\gamma_0$; then
the transformation of ${\bfH}$ corresponding to $\gamma_0$ as an
element of $\pi_1(L)$, is $z\mapsto e^A z$. Thus, the leaf $L$ is
obtained from ${\bfH}$ by identifying $z$ and $e^A z$.

Now, let us consider a typical Brownian trajectory $\gamma$ in $L$
and its lift to the universal cover $\tilde{\gamma}(t)=(x(t),y(t))$.
Note that $\tilde{\gamma}$ satisfies the following stochastic
differential equation:
$$
\left\{
\begin{array}{l}
\dot{x}= \sqrt{2} y \, dW^1_t, \\
\dot{y}= cy + \sqrt{2} y \, dW^2_t,
\end{array}
\right.
$$
where $W_t^1$ and $W_t^2$ are two independent Wiener processes. Here
the coefficient $\sqrt{2}$ comes from our definition of Brownian
motion: as we have defined it using the heat kernel, its intensivity
equals $2$ (instead of its common value~$1$). Let us make a change
of variables: let $u= \log y$, $v= x/y$. Then
\begin{multline*}
\dot{u} =(\log y)^{\cdot}= (\log \cdot)'(y) \cdot \kappa y +
\frac{1}{2}(\log \cdot)''(y) \cdot 2y^2 + \\
+(\log \cdot)'(y) \cdot \sqrt{2}y dW^1_t = (-1+\kappa) + \sqrt{2}
dW^1_t,
\end{multline*}
\begin{multline*}
\dot{v}= (\frac{x}{y})_{y}' \cdot \kappa y +
\frac{1}{2}(\frac{x}{y})_{yy}'' \cdot 2y^2 + (\frac{x}{y})_{x}'
\cdot \sqrt{2}y\, dW^1_t + (\frac{x}{y})_{y}' \cdot
\sqrt{2}y\, dW^2_t = \\
= -\kappa y \frac{x}{y^2} + \frac{2x}{y^3} y^2 + \sqrt{2} (dW^1_t +
\frac{x}{y} dW^2_t) = (2-\kappa) v + \sqrt{2}(dW^1_t + v dW^2_t).
\end{multline*}
This implies that $v$ satisfies a stochastic differential equation
$$
\dot{v}=(2-\kappa) v + \sqrt{2(1+v^2)} dW_t.
$$
Let us now make another change of variable: we denote
$\xi=f(v)=\log(v+\sqrt{1+v^2})$. Then $\xi$ satisfies the following
stochastic differential equation:
\begin{multline*}
\dot{\xi}=f'(\xi) (2-\kappa) v + \frac{1}{2}f''(\xi)
(\sqrt{2(1+v^2)})^2 + f'(\xi) \sqrt{2(1+v^2)} dW_t =
\\
= \frac{1}{\sqrt{1+v^2}} (2-\kappa) v - \frac{1}{4}
\frac{2v}{\sqrt{1+v^2}^3} (2(1+v^2)) + \sqrt{2}\,dW_t =
\\
= \frac{v}{\sqrt{1+v^2}} (1-\kappa) + \sqrt{2}\, dW_t.
\end{multline*}
For $\kappa>1$ we notice that the Brownian component of this
stochastic differential equation is constant, and the drift is
towards~$0$ with the velocity separated from zero for large $v$.
Thus, there exists a probability stationary measure for this process
on the real line. By lifting this measure to the initial cylinder
(by a product with the Lebesgue measure in $\log y$), we obtain a
stationary measure on $L$. We have constructed a harmonic measure
on~$L$.
\end{pf}

\newpage

\section{Similar foliations}\label{s:Similar}

A Laplace foliation $(\mathcal F, \Delta)$ is called
\textit{similar} if there exists a transverse continuous foliation
$\mathcal G$ of dimension $\codim(\mathcal F)$ such that the
operator $\Delta$ is invariant by $\mathcal G$; thus $\mathcal G$
preserves the metric $g$ and the drift vector field $V$.

\vspace{0.2cm}

The main goal of this part is to prove the unique ergodicity
property for a codimension $1$ similar Laplace foliation whose drift
vector field preserves the volume, and whose transverse structure is
just supposed continuous. We begin by giving examples of such
foliations.

\subsection{Some examples}

\subsubsection{Suspension}\label{suspension}
Let $(N,\overline{\Delta})$ be a compact manifold equipped with a
Laplace operator $\overline{\Delta}$, and $\rho : \pi_1(N) \rightarrow Homeo(F) $ be a
representation of its fundamental group into the group of
homeomorphisms of a compact manifold $F$. Let $\widetilde {N}$ be
the universal cover of $N$. The diagonal action of the discrete
group $\pi_1(N)$ on the product $\widetilde{N} \times F$ is
discontinuous and free. Moreover, it preserves the horizontal
foliation and the vertical fibration. Thus, the quotient $N \ltimes
_{\rho} F$ is equipped with a foliation $\mathcal F$ (quotient of
the horizontal foliation) and with a transverse fibration
$F\rightarrow M \stackrel{\pi}{\rightarrow} N$ (quotient of the
vertical fibration). Let $\Delta$ be the Laplace operator on the
leaves of $\mathcal F$ such that $(\pi_{\mathcal F})_* \Delta =
\overline{\Delta}$. By construction the foliation $(\mathcal F,
\Delta)$ is similar. Such foliations are called
\textit{suspensions}.

\subsubsection{Linear Anosov diffeomorphism}
Let $A:{\bf T}^n \rightarrow {\bf T}^n$ be a linear Anosov
diffeomorphism of the torus ${\bf T}^n = {\bf R}^n /{\bf Z}^n$.
Consider the quotient $M$ of $(0,\infty)\times {\bf T}^n$ by the
diffeomorphism $\widetilde{A} (t,x) = (2t,Ax)$: this is the fiber
bundle over the circle whose fiber is ${\bf T}^n$ and monodromy is
given by $A$. Define the foliations $\mathcal F$ and $\mathcal G$ to
be respectively the quotient of $(0,\infty)\times \mathcal F ^u$ and
of $\mathcal F^s$. We define a Laplace operator $\Delta$ on the
leaves of $\mathcal F$ of the form
\[  \Delta= \Delta_t + t^2 \frac{\partial ^2}{\partial t ^2} ,\]
where $\{\Delta_t\}_{t>0}$ is a family of linear Laplace operators
on $\mathcal F^u$ depending smoothly on $t$, and verifying the
relation
\[  \Delta_{2t} = (A|_{\mathcal F^u})_* \Delta_{t},\]
for every $t>0$. The Laplace foliations $(\mathcal F ,g)$ are
similar (the invariance of $\Delta$ by $\mathcal G$ comes from the
fact that the operators $\Delta_t$ are linear).

\begin{remark} The similar foliations by surfaces of a compact $3$-manifold
are well known~\cite{Carriere,Epstein}: in this case, the only
examples with interesting dynamics are the suspensions and the
foliations induced by a linear Anosov diffeomorphism of a $2$-torus.
However, it may exist other examples in higher dimension.
\end{remark}

\subsection{Non divergence of the leaves}\label{ss:non-divergence}

In this paragraph we prove that the leaves of a similar codimension
$1$ foliation whose drift vector field preserves the volume
$\mathrm{vol}_g$ are not diverging in a set of directions of large
measure. This has been observed by Thurston (see \cite{Calegari,Fenley} for
a topological proof).

\begin{definition} Let $\mathcal F$ be a similar foliation.
A \textit{transverse harmonic measure} on $\mathcal F$ is a family
$\{ \nu_L\}$ of measures $\nu_L$ on every $\mathcal G$-leaf $L$,
such that in a chart $B\times T$ in which $\mathcal F$ and $\mathcal
G$ are respectively the horizontal and vertical foliations, the
function
\[  p \in B \mapsto \nu (\{ p\} \times T )\in {\bf R}_+ \]
is harmonic.\end{definition}

\begin{lemma} On a codimension one similar Laplace foliation of a compact
manifold whose drift vector field preserves the leafwise volume
$\mathrm{vol}_g$, and which is minimal, there exists a transverse
harmonic measure.\label{transverse harmonic measure}\end{lemma}

\begin{pf} The result follows from the existence of a \textit{harmonic measure} for the
adjoint operator $\Delta^*$ of $\Delta$. Recall that $\Delta^*$ is
defined on every $\mathcal F$-leaf $L$ in such a way that for any
smooth functions $u,v:L\rightarrow {\bf R}$ with compact support one
has
\[  \int_L  u \Delta v d\mathrm{vol}_g = \int_L  (\Delta^* u) v d\mathrm{vol}_g.\]
An integration by parts shows that one has the following formula:
\[  \Delta^* = \Delta_g - V + \mathrm{div}_{\mathrm{vol}_g} V,\]
where $V$ is the drift vector fields of $\Delta$ (i.e. by definition
$\Delta = \Delta_g + V$). If $V$ preserves the volume
$\mathrm{vol}_g$, which means that the divergence of $V$ vanishes
identically, then the operator $\Delta ^*$ is also a Laplace
operator. In \cite{Candel}, it is proved that for such operators
there exists a harmonic measure.

Let $\mu$ be a harmonic measure on $(\mathcal F, \Delta ^*)$.
Consider a foliation box $B\times T$ in which $\mathcal F$ and
$\mathcal G$ are respectively the horizontal and vertical foliation.
Let $T'$ be an open subset of $T$. The image of the measure $\mu$ by
the projection $B\times T' \rightarrow B$ is a $\Delta^*$-harmonic
measure on $B$, because the foliation $\mathcal F$ is similar. Thus,
there is a $\Delta$-harmonic function $L_{T'} : B\rightarrow
[0,\infty)$ such that for any continuous function $f \in C^0_c(B) $,
one has $\int_{B\times T'}  f(b) d\mu(b,t) = \int
L_{T'}(b)f(b)\vol_g(b)$. Because $\mu$ is a measure, if $T_n$ are
disjoint open subsets of $T$, one has the relations $\sum _n L_{T_n}
= L_{\cup_n T_n}$; thus there exists a transverse measure on the
leaves of $\mathcal G$ such that $\nu (b\times T' ) := L_{T'} (b)$,
for every Borel subset $T'$ of $T$. This transverse measure is
$\Delta$-harmonic, by construction and the lemma is proved.\end{pf}

\begin{lemma} Let $\mathcal F$ be a codimension one similar Laplace
foliation of a compact manifold. Let $[x,y]$ be an interval in a
$\mathcal G$-orbit. There exists a uniquely defined map $I_{x,z}:
\widetilde{L_x}\times [x,y] \rightarrow M$ which maps the horizontal
Laplace foliation (given by the Laplace operator on
$\widetilde{L_x}$) on $(\mathcal F,\Delta)$, the vertical foliation
on $\mathcal G$, and the interval $\{x\}\times [x,y]$ identically on
$[x,y]$. \end{lemma}

\begin{pf} Let $\gamma:[0,1]\rightarrow \widetilde {L_x}$ be a smooth path starting at $x$. The restriction
of $I_{x,y}$ to $\gamma([0,1])\times [x,y]$ is uniquely defined, if
it exists. Let $0\leq t\leq 1$ be the supremum of those $t$ such
that $I_{x,y}$ is defined on $\gamma([0,t])\times [x,y]$. It is
clear that $t>0$, because locally we have foliation charts. Recall
that the foliation $\mathcal G$ preserves the metric $g$; thus, for
any $z\in [x,y]$ the length of the curve $I_{x,y}(\gamma([0,s]\times
z)$ equals the length of the curve $\gamma([0,s])$. This implies
that it is possible to extend the map $I_{x,y}$ on the domain
$\gamma([0,t]) \times [x,y]$, and because locally we have foliation
charts, to a domain $\gamma([0,t_+))\times [x,y] $, where $t_+ >t$.
Thus $t=1$ and the lemma is proved.\end{pf}

\vspace{0.2cm}

Recall that for any point $x\in M$, $\Gamma_x$ is the set of
continuous paths contained in the leaf through the point $x$, and
the Laplace operator $\Delta$ induces a probability measure $W_x$ on
$\Gamma_x$. Here is the main result of the paragraph, where
$J_{x,y}(p) $ is the point $I_{x,y}(p,1)$.

\begin{proposition} Let $\mathcal F$ be a Laplace similar foliation 
of codimension $1$ of a compact manifold, which is minimal,
and whose drift vector field preserves the volume $\mathrm{vol}_g$.
For every $\varepsilon >0$, there exists a constant $\delta >0$ such
that if $x$ and $y$ are two points on the same $\mathcal G$-orbit
with $d(x,y)\leq \delta$, then there exists a subset $E_{x,y}
\subset \Gamma_x$ of $W_x$-measure $1/2$ such that for any
$\gamma\in E_{x,y}$, one has
\[  \limsup _{t\rightarrow \infty}  d(\gamma(t) , J_{x,y}(\gamma(t) )\leq \varepsilon.\]
\label{nondivergence of leaves} \end{proposition}

\begin{pf} Consider a transverse harmonic measure $\nu$ constructed in 
Lemma~\ref{transverse harmonic measure}.
Because $\mathcal F$ is minimal, the measure $\nu$ restricted to
every leaf of $\mathcal G$ has full support, and is diffuse. Thus,
there exists $\delta '>0$ such that if $[x,y]$ is an interval in a
$\mathcal G$-leaf of $\nu$-measure bounded by $\delta'$, then the
distance between $x$ and $y$ is bounded by $\varepsilon$. Moreover,
there exists $\delta >0$ such that if the distance between $x$ and
$y$ is bounded by $\delta$, then the $\nu$-measure of the interval
$[x,y]$ is bounded by $\delta'/2$.

The function $p \in \widetilde {L_x} \mapsto f(\gamma)=\nu
(I_{x,y}(\{p \} \times [x,y]))\in (0,\infty)$ is harmonic. Thus by
the martingale theorem, for $w_x$-almost every $\gamma$, the limit
\[  \lim _{t\rightarrow \infty} f(\gamma(t))\]
exists, and its integral over $\Gamma_x$ is $f(x) = \nu([x,y])\leq
\delta'/2$. In particular, there is a measurable subset
$E_{x,y}\subset \Gamma_x$ of $W_x$-measure $1/2$ such that for every
$\gamma\in E_{x,y}$
\[ \lim _{t\rightarrow \infty} f(\gamma(t))\leq  \delta' . \]
For every $\gamma$ of $E_{x,y}$ we have
\[\limsup _{t\rightarrow \infty}  d(\gamma(t) , I_{x,y}(\gamma)(t,1)) \leq \varepsilon.\]
The proposition is proved.\end{pf}

\subsection{Application to unique ergodicity}
We prove that there is only one harmonic measure on a codimension
$1$ similar foliation whose drift vector field preserves the volume
$\mathrm{vol}_g$.

\begin{theorem} Let $(\mathcal F, \Delta)$ be a similar Laplace foliation of codimension $1$ of a compact manifold $M$, whose drift vector
field preserves the volume $\mathrm{vol}_g$. Then, on a minimal
subset of $\mathcal F$ is supported a unique harmonic
measure.\label{unique-ergodicity on a minimal} \end{theorem}

\begin{pf} Let $\mathcal M$ be a minimal subset of $\mathcal F$. There are three possibilities:

\begin{itemize}
\item $\mathcal M$ is a compact leaf.
\item $\mathcal M$ is transversely a Cantor set.
\item $\mathcal M = M$.
\end{itemize}
In the first case, the only harmonic measure is the unique
$\Delta$-harmonic volume on the compact leaf. The second case can be
reduced to the third one by collapsing the components of the leaf of
$\mathcal G$ outside $\mathcal M$. Thus, we suppose that $\mathcal M
= M$, i.e. $\mathcal F$ is minimal.

\begin{lemma} Let $\mathcal F$ be a similar minimal foliation of codimension $1$ of a compact manifold $M$,
and $\mu$ be an ergodic harmonic measure. Then for every $\alpha
>0$, and every continuous function $f : M\rightarrow {\bf R}$, the
following property holds. For every point $y\in M$, there exists a
measurable subset $E_y \subset \Gamma_y$ of $W_y$-measure $1/2$ such
that for every $\gamma \in E_y$:
\[\int fd\mu -\alpha \leq \liminf_{n\rightarrow \infty} B_n(f,\gamma)
\leq \limsup_{n\rightarrow \infty} B_n(f,\gamma) \leq \int fd\mu
+\alpha,\] where $B_n(f,\gamma):= \frac{1}{n} \sum _{1\leq k\leq n }
f(\gamma(k))$ are the Birkhoff sums of $f$ along the path $\gamma$.
\label{convergence of Birkhoff sums}\end{lemma}

\begin{pf} By the Birkhoff theorem for harmonic measures proved in \cite{Garnett},
there is a measurable subset $X\subset M$ of full $\mu$-measure,
saturated by $\mathcal F$, so that for every $x\in X$ and
$w_x$-almost every continuous path $\gamma \in \Gamma_x$, the
Birkhoff sums $B_n(f,\gamma)$ converge to $\int f d\mu$.

Let $\varepsilon>0$ such that if $d(x,y)\leq \varepsilon$ then
$|f(x) -f(y)|\leq \alpha$. Let $y$ be any point of $M$. There exists
a point $x\in X$ which is in the $\mathcal G$-orbit of $y$ and such
that $d(x,y)\leq \delta$, the $\delta$ being given by Lemma
\ref{nondivergence of leaves}. Let $F_x$ denote the set of element
$\gamma\in \Gamma_x$ for which the Birkhoff sums converge to $\int f
d\mu$. Because $x$ belongs to $X$ this set is of full measure.
Define $E_y:= J_{x,y} (F_x \cap E_{x,y})$. The set $E_y$ is of
$W_y$-measure $1/2$, because $F_x\cap E_{x,y}$ is of $W_x$-measure
$1/2$ and that $J_{x,y}$ sends $W_x$ on $W_y$. Let $\gamma = J_{x,y}
(\gamma')\in E_y$. Using Lemma \ref{nondivergence of leaves}, we
have
\[  \limsup _{t\rightarrow \infty} d(\gamma(t),\gamma'(t)) \leq \varepsilon.\]
Thus one gets
\[  \limsup _{n\rightarrow \infty} |B_n(f,\gamma) -B_n (f,\gamma')| \leq \alpha,\]
and the lemma follows because the Birkhoff sums of $\gamma '$
converge to $\int f d\mu$.\end{pf}

\vspace{0.2cm}

We are now able to finish the proof of the theorem. Let $\mu'$ be
another ergodic measure, and $f$ be a continuous function on $M$. We
are going to prove that $\int f d\mu ' = \int f d\mu$. Observe that
by ergodicity of $\mu'$ there exists a point $y$ on $M$ such that
for $w_y$-almost every $\gamma\in \Gamma_y$, the Birkhoff sums
$B_n(f,\gamma)$ converge to $\int f d\mu'$. Denote by $F_y$ the set
of such paths $\gamma$ which is of full $w_y$-measure. We apply
Lemma \ref{convergence of Birkhoff sums} to the point $y$: for every
$\gamma\in E_y$, the Birkhoff sums $B_n(f,\gamma)$ are tending to
$\int f d\mu$ with an error of $\alpha$. Because the measure of
$E_y$ is positive, it intersects $F_y$. Thus one gets $|\int f d\mu'
-\int f d\mu| \leq \alpha$. Because $\alpha$ is arbitrary, there is
only one ergodic harmonic measure, and thus only one harmonic
measure. The theorem is proved.\end{pf}

\begin{proposition} Let $(\mathcal F, \Delta_g)$ be a similar Laplace foliation of codimension $1$ of
a compact manifold $M$, where $\Delta_g$ is the Laplacian of a
riemannian metric. Then every ergodic harmonic measure is supported
on a minimal subset.\end{proposition}

\begin{pf} Let $\mu$ be an ergodic harmonic measure on $(\mathcal F, \Delta_g)$, and $\mathcal M$ be a minimal closed subset
contained in the support of $\mu$. We are going to prove that
$\mathcal M$ is exactly the support of $\mu$. It is clear that one
can suppose that $\mu$ does not charge any leaf.

Suppose that the foliation is oriented. Because the operator
$\Delta_g$ is symmetric, the measure $\mu$ induces a transverse
harmonic measure which is $\Delta_g$-harmonic. From every point $x$
of $M$, the positive $\mathcal G$-orbit of $x$ intersects $\mathcal
M$ in a first time in a point $y$. Consider the function $f(x) =
\nu([x,y))$. This is a continuous function, because $\mu$ does not
charge any leaf, and it is harmonic on every leaf. By Garnett lemma
\cite{Garnett}, this function has to be constant on $\mu$-almost
every leaf, thus by continuity of $f$, $f$ is constant on the
support of $\mu$. But on the minimal $\mathcal M$, $f$ vanishes, so
that the restriction of $f$ to the support of $\mu$ is identically
$0$. Thus the support of $\mu$ has to be reduced to $\mathcal M$.
The proposition is proved.\end{pf}

\begin{example}\label{ex:lamination}
Theorem~\ref{unique-ergodicity on a minimal} seems to be false
when the drift does not preserve the volume $\mathrm{vol}_g$. We
give an example of a similar Laplace \textit{lamination} of a
compact space which is minimal, transversely conformal, and not
uniquely ergodic.

A lamination of a compact space $X$ is an atlas of
homeomorphisms from open sets of $X$ to the product of an euclidian
ball by a topological set, in such a way that the changes of
coordinates preserve the local fibration by balls, and the
diffeomorphisms from a piece of ball to another depends
continuously of the transverse parameter in the smooth topology. The
definition of a Laplace operator on a lamination is exactly the same
as in the foliation case.

We are going to describe an example of a similar Laplace lamination
of a compact space, which has been constructed by Sullivan~\cite{Su3}.
Let ${\bf H}$ be the upper-half plane, whose hyperbolic metric is
expressed by
\[  g = \frac{dx^2 + dy^2}{y^2},\]
in the coordinates $z=x+iy$ of ${\bf H}$. Consider the unit vector
fields that points on the direction $\infty$ of $\partial {\bf H}$.
In the $x,y$ coordinates, it is expressed as
\[  V = y\frac{\partial}{\partial y}.\]
The direct isometries of $({\bf H},g,V)$ are the maps of the form
$z\mapsto a z + b $ where $a$ is a positive number, and $b$ is a
real number. For any real number $\kappa$, these transformations
preserve the Laplace operator
\[  \Delta_{\kappa} = \Delta_g + \kappa V.\]

Let $A({\bf Z}[1/2])$ be the group of affine transformations
$x\mapsto ax + b$, where $a$ is a power of $2$, and $b$ is a dyadic
integer of the form $p/2^n$, $p$ and $n$ being integers. The group
$A({\bf Z}[1/2])$ acts naturally on the product ${\bf H} \times {\bf
Q}_2$ of the upper half plane by the field of $2$-adic numbers, the
action preserving the natural structure of the horizontal Laplace
lamination $({\bf H}\times {\bf Q}_2,\Delta_{\kappa})$. The action
is discrete, without fixed point, and the quotient
$(X,\Delta_{\kappa})$ is a similar Laplace lamination of a compact space.

When $\kappa > 1$ the Laplace laminations
$(X,\Delta_{\kappa})$ carry harmonic measures that charge any
cylinder leaves. This can be seen by the same arguments as those
given in \ref{s:counter-example}.
These examples are not codimension $1$ foliations,
but they share with codimension $1$ foliation the property of being
transversely conformal, which is the only property used in the proof
of Theorem \ref{unique-ergodicity on a minimal}.

It seems to us that the hypothesis $div_{\mathrm{vol}_g}
V=0$ is too strong. For instance, we conjecture that for any Laplace
operator on the base of a suspension, the conclusion of
\ref{unique-ergodicity on a minimal} holds.

There are analog examples of similar Laplace laminations of a
compact space, associated with tilings of the hyperbolic plane. They
have been studied by Petite \cite{Petite}; many of them are minimal
but not uniquely ergodic even for a symmetric operator. The lack of
unique ergodicity is due to the fact that they do not carry a
transversely invariant ``conformal" structure.
\end{example}

\newpage
\section{Appendix: Technical proofs}\label{s:TAP}

\subsection{}

\begin{pfof}{Proposition~\ref{p:discrcomp}}
It is a well-known fact (see~\cite{CLY,Malliavin}), that for a
Brownian motion, the probability of making large steps decreases
very fast:
\begin{equation}
    \exists C_1,d_0: \quad \forall p\in M, \,
    \forall d>d_0 \quad
    W_p(\dist(\widetilde{\gamma}(\delta), p)>d)\le e^{-C_1 d}.
    \label{Wexp}
\end{equation}
Thus, we may choose a random variable $\xi$, such that $\xi\ge 0$,
$E\xi<\infty$ and
\begin{equation}\label{Wbound}
\forall p\in M, \forall d>0 \quad
W_x(\dist(\widetilde{\gamma}(\delta), p)>d) \le P(\xi>d).
\end{equation}
Let us choose for every $x\in\mcM$, a function
$\chi_x:\Gamma_x\to\bbR$, such that $\chi_x$ depends only on
$\gamma|_{[0,\delta]}$, has the same distribution as $\xi$ and such
that $\chi_x(\gamma)\ge k_{1}(\gamma)$ for every
$\gamma\in\Gamma_x$. Recall, that $\sigma:\Gamma\to\Gamma$ stays for
the map, erasing the first step of the $\delta$-discretization:
$\sigma(\gamma)(t)=\gamma(t+\delta)$. Then,
$$
\forall j,\gamma \quad k_{j}(\gamma)=k_{1}(\sigma^{j-1}(\gamma)).
$$
Let us denote
$\xi_j(\gamma)=\chi_{\gamma((j-1)\delta)}(\sigma^{j-1}(\gamma))$.
Then,
\begin{multline*}
\frac{k_{1}(\gamma)+\dots+k_{n}(\gamma)}{n} = \frac{k_{1}(\gamma)+
k_{1}(\sigma(\gamma))+ \dots+ k_{1}(\sigma^{n-1}\gamma)}{n} \le
\\
\le \frac{\chi_{\gamma(0)}(\gamma)+
\chi_{\gamma(\delta)}(\sigma(\gamma))+ \dots+
\chi_{\gamma((n-1)\delta)}(\sigma^{n-1}\gamma)}{n} =
\\
= \frac{\xi_{1}(\gamma) + \dots + \xi_{n}(\gamma)}{n}
\end{multline*}
Note, that all the $\xi_j$ are distributed identically with $\xi$.
Moreover, from the Markovian property, the conditional distribution
of the variable $\xi_{n+1}$ with respect to every condition
$\gamma|_{[0,n\delta]}=\overline{\gamma}$ coincides with the
distribution of~$\xi$. On the contrary, the variables
$\xi_1,\dots,\xi_n$ are determined by such a condition. Thus, the
variable $\xi_{n+1}$ is independent from $\xi_1,\dots,\xi_n$. As $n$
is arbitrary, all the variables $\xi_1,\dots,\xi_n,\dots$ are
independent and identically distributed.

Now, let us take $c=E\xi +1$. We are going to show that for almost
every $\gamma\in\Gamma$ we have $\limsup\limits_{n\to\infty} K_n
(\gamma)/n <c$. From the Large Numbers Law, we have:
$$
\limsup\limits_{n\to\infty}
\frac{\xi_1(\gamma)+\dots+\xi_n(\gamma)}{n} = E\xi < c
$$
$W_x$--almost surely. But
$$
\limsup\limits_{n\to\infty} \frac{k_1(\gamma)+\dots+k_n(\gamma)}{n}
\le \limsup\limits_{n\to\infty}
\frac{\xi_1(\gamma)+\dots+\xi_n(\gamma)}{n},
$$
and thus
$$
\limsup\limits_{n\to\infty} \frac{k_1(\gamma)+\dots+k_n(\gamma)}{n}
< c
$$
$W_x$-almost surely. This concludes the proof of the proposition.
\end{pfof}

\subsection{Proof of Lemma~\ref{l:positive}}

First, we are going to prove the lemma in the particular case of a
codimension one foliation~$\mcF$, using the transversal
one-dimensional foliation~$\mcG$, described in
Section~\ref{s:Negative}.

In order to prove the lemma, we shall study the behaviour of heat
kernels for a time $t=\delta$ fixed on different but close enough
leaves. We are going to use the following idea: major parts of the
heat distribution measures on these two leaves are similar (i.e. the
density of $\mcG$-holonomy image of one with respect to another is
close to~$1$). Moreover, the infinite product of these densities
converges along most Brownian paths, because along these paths
the leaves approach exponentially. Thus, the measures
$W_x^{\delta}=F_* W_x|_{E_x}$ and $(F \circ \Phi_{\bar{x},x})_*
W_{\bar{x}}|_{E_{\bar{x}}}$ are absolutely continuous with respect
to each other, and the density on the major part of trajectories is
close to~$1$. In particular, the total measures $W_x(E_x)$ and
$W_{\bar{x}}(E_{\bar{x}})$ are close to each other.


\begin{proposition}\label{p:hkb}
Let $\bar{x}\in\mcG_x$, $\dist_{\mcG}(x,\bar{x})=\theta$. Consider
the measures $\nu_x$ and $\nu_{\bar{x}}$, where the measure
$\nu_z=p(z,\cdot;\delta)\, d\vol_g$ on $\mcF_z$ gives the
distribution of Brownian motion at the time~$\delta$. Also, let
$\varepsilon_2>0$ and $R$ be chosen. Then, there exists a set
$S=S(R,\varepsilon_2,\theta,x)\subset B_{R}^{\mcF}(x)$, such that
\begin{itemize}
\item $\nu_x(S)>1-\varepsilon_1$,
\item $\left.\frac{d\nu_x}{d (\Phi_{x,\bar{x}}^* \nu_{\bar{x}})}\right|_{S}\in
[1-\varepsilon_2,1+\varepsilon_2]$,
\end{itemize}
where $\varepsilon_1= \frac{1}{\varepsilon_2} \left( G_1' e^{G_2' R}
\theta+ e^{- G_3 R^2} \right)$, and $G_1', G_2', G_3$ are geometric
constants (depending only on the foliation $\mcF$).
\end{proposition}
\begin{pf}
First, let us equip the leaf $\mcF_x$ with another Riemannian metric
$g'$, coinciding with $\Phi_* g|_{\mcF_{\bar{x}}}$ inside
$B_{R}^{\mcF}(x)$ and with $g|_{\mcF_x}$ outside $B_{2R}^{\mcF}(x)$;
in the annulus left, the metrics $g'$ is defined using cut-off
function.

Note, that the distance (in $C^{2d}$-topology, where $d$ is the
dimension of the leaves of~$\mcF$) between $g$ and $g'$ is at
most~$G_0 \theta e^{G_1 \cdot 2R},$ where $G_1$ is a constant,
giving the maximum deviation of leaves of~$\mcF$, and $G_0$ gives
the maximum of derivative of $g$ along~$\mcG$.

Also, we recall, that the heat kernel can be constructed explicitly
as a series of convolutions. This procedure is described in the
books of Candel~\cite{Candel1} and Chavel~\cite{Chavel}. In a few
words, a function $L$, which ``almost satisfies'' the heat equation,
is explicitly constructed, and then the real heat kernel is obtained
as a sum of $L$ and a series of convolutions. These series are
converging uniformly for every fixed moment of time $t$, and the
dependence of the metrics is smooth (due to the explicit nature of
the construction). Thus, the distance between the heat kernels for
the metrics $g$ and $g'$ at the time $t=\delta$ can be bounded by
the product of a constant $G_2$ (depending only on the geometry of
foliation $\mcF$ and of the moment $\delta$) and of the distance
between $g$ and~$g'$. This distance is at most $e^{2 G_1 R} \theta$;
hence, the difference between these kernels is bounded by $ G_2 e^{2
G_1 R}\theta.$

Now, due to the upper bounds for the heat
kernel~\cite{CLY,Malliavin}, the set of Brownian trajectories on the
interval of time $[0,\delta]$, starting at $x$ and exiting from the
ball $B_{R}^{\mcF}(x)$ at some intermediate moment, has the measure
at most $e^{-G_3 R^2}$. We notice that these trajectories for the
metric $g'$ and for the metric $\Phi_{x,\bar{x}}^{*}
g|_{\mcF_{\bar{x}}}$ are the same (they do not pass through the
points where these metrics do not coincide). Thus, the parts of heat
kernel at time $\delta$, coming from these trajectories, are the
same for these two metrics.

Let us define the set $S$ in a following way: $z\in S$, if
\begin{itemize}
\item The density $p_{g'}(x,z;\delta)$ is at least
$3G_2 e^{2 G_1 R} \theta / \varepsilon_2$
\item At least $1-\varepsilon_2/3$ of this density comes from the
trajectories staying inside $B_R^{\mcF}$ (in particular, $z\in
B_R^{\mcF}(x)$).
\end{itemize}
Then, the second condition for $S$ (quotient of densities) is
satisfied automatically: the maximum possible change of density at a
point of $S$ is the sum of changes while passing from $g$ to $g'$
(at most $\varepsilon_2/3$ part) and from $g'$ to
$g|_{\mcF_{\bar{x}}}$ (at most $\varepsilon_2/3$ part due to common
set of trajectories).

Now, let us estimate $\nu(\mcF_x\setminus S)$, thus verifying the
first condition. The points of this complementary can be of two
types: either points of $B_R^{\mcF}$ with too small value of density
(we note this set $X_1$), or with too big part of this density
coming from trajectories, exiting the ball~$B_R^{\mcF}(x)$ (we note
this set $X_2$).

The first part is estimated as
\begin{multline*}
\nu(X_1)= \int_{B_R^{\mcF}(x)} p(x,z;\delta) d\vol_g \le 3G_2 e^{2
G_1 R} \theta / \varepsilon_2 \cdot \vol_g(B_R^{\mcF}(x)) \le
\\
\le (3G_2 e^{2 G_1 R} \theta / \varepsilon_2) \cdot e^{G_4 R},
\end{multline*}
where $G_4$ is the constant, bounding the growth of the leaves
of~$\mcF$.

The second part is estimated as follows: denote by $\rho(z)$ the
part of the density $p(x,z;\delta)$, coming from the trajectories
exiting from the ball $B_R^{\mcF}(x)$. Then,
\begin{multline*}
\nu(X_2)= \int_{\{z:\, \rho(z)/p(x,z;\delta)>\varepsilon_2/3\}}
p(x,z;\delta) \, d\vol_g(z) \le
\\
\le \int_{\{z:\, \rho(z)/p(x,z;\delta)>\varepsilon_2/3\}}
\frac{3}{\varepsilon_2} \cdot \rho(z) \, d\vol_g(z) \le
\\
\le \frac{3}{\varepsilon_2} \int_{\mcF_x} \rho(z) \, d\vol_g(z) <
\frac{3}{\varepsilon_2} \cdot e^{-G_3 R^2}
\end{multline*}
(the last inequality comes from the upper bound for the probability
of all the set of trajectories, leaving the ball of radius $R$
at some moment between 0 and~$\delta$).

Finally, we obtain
\begin{multline*}
\nu(\mcF_x\setminus S)= \nu(X_1\cup X_2) < \frac{ 3G_2 e^{(2 G_1 +
G_4) R}}{\varepsilon_2} \cdot \theta + \frac{3}{\varepsilon_2} \cdot
e^{-G_3 R^2} =
\\
=\frac{1}{\varepsilon_2} \left( G_1' e^{G_2' R} \theta+ e^{- G_3
R^2} \right),
\end{multline*}
where $G_1'=3G_2$, $G_2'=2G_1+G_4$.

The first condition on~$S$ is satisfied.
\end{pf}

\vspace{1cm}

For any $R,\theta>0$ let us denote
$$
\Psi(R,\theta)=\sqrt{G_1' e^{G_2' R} \theta+ e^{- G_3 R^2}}.
$$
Also, let $r(\theta)=(\log \frac{1}{2\theta G_1'})/(2G_2')$. Then
$G_1' e^{G_2' r(\theta)}\theta= \frac{1}{2} \sqrt{\theta}$, and for
all $\theta$ sufficiently small $e^{-G_3 r(\theta)^2}< \frac{1}{2}
\sqrt{\theta}$. Thus, for all $\theta$ sufficiently small
$$
\Psi(\theta):=\Psi(r(\theta),\theta)<\sqrt{\theta}.
$$
Now, denote
$$
S(\theta,x)=S(r(\theta),\Psi(\theta),\theta,x).
$$
For this set, the conclusions of Proposition~\ref{p:hkb} is
satisfied with $\varepsilon_1=\varepsilon_2<\sqrt{\theta}$. Let us
choose a small transversal interval $J\subset I$, $J\ni x$, and
consider the subset of $E_x$, defined as
\begin{equation}
E_x'=\left\{ \gamma\in E_x | \, \forall n\ge 0 \quad x_{n+1}\in
S(\theta_n, x_n) \right\},
\end{equation}
where $x_n=\gamma(n\delta)$ is the discretization sequence,
corresponding to~$\gamma$, and $\theta_n=|h_{\gamma|_{[0,n\delta]}}
(J)|$ is the sequence of the corresponding transverse distances
(exponentially decreasing due to the nature of $E_x$).

\begin{lemma}\label{l:ac}
    For $\bar{x}\in J$, the images of the measure $W_x|_{E_x'}$ under $F$ and of
    the measure $W_{\bar{x}}|_{E_{\bar{x}}'}$ under $F\circ \Phi_{\bar{x},x}$ are
    absolutely continuous with respect to each other. The density
    can be made arbitrary close to~1 by choice of sufficiently small
    intervals~$I$ and $J$. Moreover, by such a choice,
    the differences $W_x(E_x)-W_x(E_x')$ and
    $W_{\bar{x}}(E_{\bar{x}})-W_{\bar{x}}(E_{\bar{x}}')$ can be made
    arbitrarily small.
\end{lemma}

\begin{pf}
Let us consider the projection map
$$
\pi_n:(\widetilde{\mcF}_{x})^{\infty}\to
(\widetilde{\mcF}_{x})^{n+1}, \quad
\pi_n(\{x_j\}_{j=0}^{\infty})=\{x_j\}_{j=0}^{n}
$$
and its composition with $F$, which we denote
$F_n:E_x'\to(\widetilde{\mcF}_x)^{n+1}$,
\begin{equation*}
    F_n(\gamma)=\{\widetilde{\gamma}(j\delta)\}_{j=0}^{n}.
\end{equation*}

Let us denote
$$
\mu_1=F_* W_x|_{E_x'}, \quad \mu_2=(F\circ \Phi_{\bar{x},x})_*
W_{\bar{x}}|_{E_{\bar{x}}'},
$$
$$
\mu_1^n=(\pi_n)_{*} \mu_1 = (F_n)_{*} W_x|_{E_x'}, \quad
\mu_2^n=(\pi_n)_{*} \mu_2 = (F_n \circ \Phi_{\bar{x},x})_{*}
W_{\bar{x}}|_{E_{\bar{x}}'}.
$$

\begin{proposition}
    The measures $\mu_1^n$ and $\mu_2^n$ are absolutely continuous
with respect to each other, and
\begin{equation*}
    \left.\frac{d\mu_2^n}{d\mu_1^n}\right|_{(x_j)}=
    \rho_n ((x_j)) =
    \prod\limits_{j=0}^{n-1}
    \frac{d\nu_{x_j}}{d(\Phi_{x_j,\bar{x}_j}^* \nu_{\bar{x}_j})},
\end{equation*}
\end{proposition}

To prove Lemma~\ref{l:ac}, it suffices to show that for
$\mu_1$--almost every point $(x_j)\in (\widetilde{\mcF}_x)^{\bbN}$
the sequence $\rho_n((x_j))$ converges to some number between $0$
and~$\infty$. It is equivalent to the convergence of the infinite
product
\begin{equation}\label{eq:prod}
\prod\limits_{j=1}^{\infty}
    \frac{d\nu_{x_j}}{d(\Phi_{x_j,\bar{x}_j}^* \nu_{\bar{x}_j})}=
    \lim_{n\to\infty} \rho_n((x_j))
\end{equation}
(again, for $\mu_1$--almost every $(x_j)\in
(\widetilde{\mcF}_x)^{\bbN}$).

Now, for every $(x_j)$ in the image $F(E_x')$, the logarithm of the
product~\eqref{eq:prod} can be estimated as
\begin{multline}\label{eq:density}
\left| \log \prod\limits_{j=1}^{\infty}
    \frac{d\nu_{x_j}}{d(\Phi_{x_j,\bar{x}_j}^* \nu_{\bar{x}_j})}
    \right| \le \sum_{j=0}^{\infty}
    \left| \log \frac{d\nu_{x_j}}{d(\Phi_{x_j,\bar{x}_j}^* \nu_{\bar{x}_j})}
    \right|
    \\
    \le \sum_{j=0}^{\infty}
    \left| 2\sqrt{\theta_j} \right| \le \sum_{j=0}^{\infty}
    \left| 2\sqrt{C e^{-\alpha j} |J|} \right| \le C_3 \sqrt{|J|}
    \sum_{j=0}^{\infty} e^{-\alpha j/2} = C_4 \sqrt{|J|}.
\end{multline}
Here we used the definition of $E_x'$ to bound the density
$\frac{d\nu_{x_j}}{d(\Phi_{x_j,\bar{x}_j}^* \nu_{\bar{x}_j})}$, and
then again to estimate $\theta_n$.

We have estimated the density; moreover, for $J$ sufficiently small
this density (due to~\eqref{eq:density}) can be made arbitrary close
to~$1$.

Now, let us estimate the difference $W_z(E_z)-W_z(E_z')$, where
$z\in J$. Denote
$$
E_{z,n}=\left\{ \gamma\in E_z | \, \forall j, 0\le j<n \quad
x_{n+1}\in S(\theta_n, x_n) \right\},
$$
$$
\widetilde{C}_{z,n}= F_n(E_{z,n}),
$$
and let $\widetilde{\mu}_n$ be a measure on
$(\widetilde{\mcF}_z)^{n+1}$, defined as the discretization image
of~$W_z|_{E_{z,n}}$. Also, consider projection maps
$\widetilde{\pi}_n: \widetilde{C}_{z,n+1}\to \widetilde{C}_{z,n}$.
Then, $(\widetilde{\pi}_n)_* \widetilde{\mu}_{n+1}$ is absolutely
continuous with respect to $\widetilde{\mu}_n$, and the density is
equal to
$$
\widetilde{\rho}_n ((x_j))=\nu_{x_n} (S(\theta_n,x_n)) \ge
1-2\sqrt{\theta_n} \ge 1- 2\sqrt{Ce^{-\alpha n} |J|}.
$$
Thus,
\begin{multline*}
W_z(E_{n,z}\setminus E_{n+1,z}) = \int_{\widetilde{C}_{n,z}}
(1-\widetilde{\rho}_n)((x_j)) \, d \widetilde{\mu}_{n}
((x_j)_{j=0}^n) \le
\\
\le \int\limits_{\widetilde{C}_{n,z}} 2\sqrt{Ce^{-\alpha n} |J|}\,
d\widetilde{\mu}_n ((x_j)_{j=0}^n) \le C_3 e^{-\alpha n/2}
\sqrt{|J|}.
\end{multline*}
Now, recall that $E_z'=\bigcap_{n=1}^{\infty} E_{z,n}$, hence,
\begin{multline}\label{eq:measure}
W_z(E_z)-W_z(E_z')=\sum_{n=0}^{\infty} W_z(E_{z,n}\setminus
E_{z,n+1}) \le
\\
\le \sum_{n=0}^{\infty} C_3 e^{-\alpha n/2} \sqrt{|J|} \le C_4
\sqrt{|J|}.
\end{multline}
The difference $W_z(E_z)- W_z(E_z')$ tends to 0 as $|J|$ tends to
$0$.
\end{pf}

\begin{pfof}{Lemma~\ref{l:positive}}
\subsubsection{``Positive measure'' part}
    First, let us prove the ``positive measure'' part. Namely,
    estimate the difference
    $|W_x(E_x)-W_{\bar{x}}(E_{\bar{x}})|$ if $\bar{x}\in J\subset
    I$:
    \begin{multline*}
    |W_x(E_x)-W_{\bar{x}}(E_{\bar{x}})| \le  |W_x(E_x)-W_x(E_x')| +
    \\
    + |W_x(E_x')-W_{\bar{x}}(E_{\bar{x}}')| +
    |W_{\bar{x}}(E_{\bar{x}}')-W_{\bar{x}}(E_{\bar{x}})|.
    \end{multline*}
    All the three differences can be estimated using Lemma~\ref{l:ac}.
    Thus, choosing any $p_1<W_x(E_x)$, we can find a sufficiently small
    transversal interval $J$, such that for any $\bar{x}\in J$ we have
    $W_{\bar{x}}(E_{\bar{x}})\ge p_1$.

    Now, suppose that
    $y\in \widetilde{\mcF}_{\bar{x}}$.
    Then, the measures $\nu_{y}$ and $\nu_{\bar{x}}$ are
    absolutely continuous with respect to each other, and the
    density on the major (with respect to these measures) part of
    $\widetilde{\mcF}_{\bar{x}}$ is close to~1. Let us choose
    $\varepsilon_3>0$ and denote the set
    $$
    \widetilde{S}(\bar{x},y,\varepsilon_3)=\left\{z\in\widetilde{\mcF}_{\bar{x}}
    \mid \frac{d\nu_{\bar{x}}}{d\nu_{y}}(z)\in
    [\frac{1}{1+\varepsilon_3},1+\varepsilon_3] \right\}.
    $$
    Then,
    $$
    \forall \varepsilon_3>0 \, \exists r>0 : \, \forall
    y\in\widetilde{\mcF}_{\bar{x}}, d(y,\bar{x})<r \quad
    \nu_{\bar{x}}(\tilde{S})>1-\varepsilon_3,
    \nu_{y}(\tilde{S})>1-\varepsilon_3.
    $$
    Note that due to the Markovian property (the conditional
    distribution of $\sigma(\gamma)$ with respect to every condition
    $\sigma(\gamma)(0)=z'$ coincides with $W_{z'}$), we have
    $$
        W_{z}(E_z)=\int_{\widetilde{\mcF}_z} W_{z'} (E_{z'}^1) \,
        d\nu_z(z'),
    $$
    where
    \begin{equation*}
    E_{z'}^1=\left\{\gamma \mid \forall n \quad
    \left| h_{\gamma|_{[0,n\delta]}}(h_{z,z'}(I)) \right| < C_0 e^{-\alpha \cdot
    (n+1)\delta} |I| \right\}
    \end{equation*}
    (this is the definition of $E_z$, rewritten in terms of the shift $\sigma(\gamma)$).
    Thus, for $d(y,\bar{x})<r$ we have
    \begin{equation}\label{eq:dstep}
        W_{y}(E_y)=W_{y}(E_y\cap\{x_1\notin \widetilde{S}\})
        + \int_{\widetilde{S}} W_z (E_z^1) \,
        d\nu_{y}(z),
    \end{equation}
    \begin{multline}\label{eq:s-change}
        W_{\bar{x}}(E_{\bar{x}})=W_{\bar{x}}(E_{\bar{x}}\cap\{x_1\notin
        \widetilde{S}\})
        + \int_{\widetilde{S}} W_z (E_z^1) \,
        d\nu_{\bar{x}}(z)=
        \\
        = W_{\bar{x}}(E_{\bar{x}}\cap\{x_1\notin \widetilde{S}\})
        + \int_{\widetilde{S}} W_z (E_z^1)
        \frac{d\nu_{\bar{x}}}{d\nu_y}(z) \, d\nu_y(z).
    \end{multline}
    The first summands in the right hand sides of~\eqref{eq:dstep}
    and~\eqref{eq:s-change} are no greater than $\varepsilon_3$, and the
    quotient of the second summands is bounded by $1+\varepsilon_3$.
    Thus, for $y$ and $\bar{x}$ being sufficiently close to each
    other, the probabilities $W_{\bar{x}}(E_{\bar{x}})$ and
    $W_y(E_y)$ are also close to each other. In particular, for
    every $p_0<p_1$ we can find $r>0$, such that for $U$
    being the union of $r$-leafwise-neighborhoods of points of
    $J$, we have
    $$
    \forall y\in U \quad W_y(E_y)\ge p_0.
    $$
    It completes the proof of this part of the lemma.

\subsubsection{Proof of the ``distributions'' part}

  Note that due to the arguments already used we may suppose that
  $y\in\mcG_x$: for $y$ and $y'$ on the same $\mcF$-leaf, the
  measures $\nu_{y'}$ and $\nu_y$ are absolutely continuous with
  respect to each other, and hence any tail-type property holds (or
  does not hold) simultaneously for typical trajectories in $\Gamma_y$ and
  $\Gamma_{y'}$. Hence, if $y$ does not belong to
  $\mcG_x$, we may replace it by $y'$, which is an intersection
  point of $\mcF_y$ and $\mcG_x$.

  Recall that almost every trajectory $\gamma\in\Gamma_x$
  (due to the choice of $x$) is
  distributed with respect to~$\mu$.
  Thus, for almost every path $\gamma\in\Gamma_x$ we have
  \begin{equation}\label{eq:path_distribution}
  \lim\limits_{T\to\infty} \frac{1}{T} \gamma_*\leb_{[0,T]} =\mu.
  \end{equation}

We know that a trajectory of $E_{y}$ approaches a trajectory of
$E_x$. Unfortunately, we can not claim that the map giving the
trajectory of $E_x$ by a trajectory of $E_y$ is absolutely
continuous (or, what is the same, maps typical trajectories to
typical ones): we have this statement only for discretizations of
trajectories.

Thus, we have to prove the distributions property using
discretizations behaviour. The following arguments are a technical
realization of this idea.

  Rewrite~\eqref{eq:path_distribution} in the terms of discretization.
  Let a continuous function $\varphi$ on $M$ be chosen. Then, for $W_x$-almost
  every $\gamma\in\Gamma_x$,
  \begin{equation}\label{eq:intlimit}
  \left(\frac{1}{T} \gamma_*\leb_{[0,T]}\right)(\varphi) = \frac{1}{T} \int_0^T
  \varphi(\gamma(t))\, dt.
  \end{equation}
  It is clear that we can restrict to the moments of time of the form $T=n\delta$;
  for such $T$,
  \begin{multline}\label{eq:int}
  \frac{1}{n\delta} \int_0^{n\delta} \varphi(\gamma(t))\, dt = \frac{1}{n\delta}
  \sum\limits_{j=0}^{n-1} \int_{j\delta}^{(j+1)\delta}
  \varphi(\gamma(t))\, dt =
  \\
  = \frac{1}{n}\sum\limits_{j=0}^{n-1} \frac{1}{\delta}
  \int_{j\delta}^{(j+1)\delta} \varphi(\gamma(t))\, dt.
  \end{multline}

  Let~$z\in M$ be some point, and let us rewrite~\eqref{eq:int} for
  a typical trajectory $\gamma\in\Gamma_z$. Namely, we
  divide this sum into discrete averaging and the rest term:
  \begin{multline}\label{eq:int2}
  \frac{1}{n\delta} \int_0^{n\delta} \varphi(\gamma(t))\, dt =
  \frac{1}{n}\sum\limits_{j=0}^{n-1} \varphi(\gamma(j\delta)) +
  \\
  +  \frac{1}{n}\sum\limits_{j=0}^{n-1} \frac{1}{\delta}
  \int_{j\delta}^{(j+1)\delta} (\varphi(\gamma(t)) -\varphi(\gamma(j\delta)))
  \, dt.
  \end{multline}
  We estimate the second term in the right hand side
  of~\eqref{eq:int},
  decomposing the sum in two parts, the one
  corresponding to $j$ with $\diam(\gamma([j\delta,(j+1)\delta]))<r$
  and the one with $\diam(\gamma([j\delta,(j+1)\delta]))\ge r$.
  \begin{multline}\label{eq:int3}
  \left| \frac{1}{n}\sum\limits_{j=0}^{n-1} \frac{1}{\delta}
  \int_{j\delta}^{(j+1)\delta} (\varphi(\gamma(t)) -\varphi(\gamma(j\delta)))
  \, dt. \right| \le
  \\
  \le \frac{1}{n} \sum\limits_{j<n,
  \diam(\gamma([j\delta, (j+1)\delta]))<r}
  \frac{1}{\delta} \int_{j\delta''}^{(j+1)\delta''}
  \left| \varphi(\gamma(j\delta)) - \varphi(\gamma(t))\, dt \right| +
  \\
  + \frac{1}{n} \sum\limits_{j<n,
  \diam(\gamma([j\delta, (j+1)\delta]))\ge r}
  \frac{1}{\delta} \int_{j\delta''}^{(j+1)\delta''}
  \left| \varphi(\gamma(j\delta)) - \varphi(\gamma(t))\, dt \right| \le
  \\
  \le \omega_\varphi(r) + 2 \frac{\#\{j: \diam\ge r\}}{n} \sup_{\mcM}
  |\varphi|,
  \end{multline}
  where $\omega_\varphi$ is the modulus of continuity of the
  function~$\varphi$.

  Let us pass in~\eqref{eq:int3} to the upper limit:
  \begin{multline}\label{eq:limsupestimate}
  \limsup_{n\to\infty} \left| \frac{1}{n}\sum\limits_{j=0}^{n-1} \frac{1}{\delta}
  \int_{n\delta}^{(j+1)\delta} (\varphi(\gamma(t)) -\varphi(\gamma(j\delta)))
  \, dt. \right| \le
  \\
  \le \omega_\varphi(r)+ 2 \sup_{\mcM} |\varphi| \cdot
  \limsup_{j\to\infty} \frac{\#\{j<n \mid
  \diam(\gamma[j\delta, (j+1)\delta])\ge r\}}{n}.
  \end{multline}

  For a function~$\varphi$ chosen, the first summand can
  be made arbitrarily small by a choice of $r$
  due to the continuity of~$\varphi$. For every $r>0$ chosen, the
  second summand can be made arbitrarily small for almost all trajectories
  by a choice of sufficiently small $\delta>0$ (which is uniform in $z\in M$)
  due to the same arguments
  as the ones used in the proof of Proposition~\ref{p:discrcomp}.

  Hence, for every $\varepsilon>0$ we can find $r$ and then $\delta$ small enough,
  such that for almost all $\gamma\in \Gamma_z$,
  $$
  \limsup_{n\to\infty} \left| \frac{1}{n}\sum\limits_{j=0}^{n-1} \frac{1}{\delta}
  \int_{j\delta}^{(j+1)\delta} (\varphi(\gamma(t)) -\varphi(\gamma(j\delta)))
  \, dt \right|<\varepsilon.
  $$

  The rest term in~\eqref{eq:int2} is estimated, and
  taking it together with~\eqref{eq:intlimit},
  for a discretization $(x_j)=F^{\delta}(\gamma)$ of a typical path
  $\gamma\in\Gamma_x$, we have
  $$
  \limsup_{n\to\infty} \left| \frac{1}{n}\sum\limits_{j=0}^{n-1}
  \varphi(x_j) - \int_M \varphi\, d\mu
  \right|<\varepsilon.
  $$

  Now, for a $W_{\bar{x}}$-typical path $\gamma$ from $E_{y}'$ let us estimate the
  difference
  \begin{equation}\label{eq:integral_estimate1}
  \limsup_{n\to\infty} \left| \frac{1}{n\delta} \int_0^{n\delta}
  \varphi(\gamma(t)) \, dt - \int_{\mcM} \varphi \, d\mu \right|.
  \end{equation}
  We have:
  \begin{multline*}
  \limsup_{n\to\infty} \left| \frac{1}{n\delta} \int_0^{n\delta}
  \varphi(\gamma(t)) \, dt - \int_{\mcM} \varphi \, d\mu \right| \le
  \\
  \le \limsup_{n\to\infty} \left| \frac{1}{n\delta} \int_0^{n\delta}
  \varphi(\gamma(t)) \, dt -  \frac{1}{n}\sum_{j=0}^{n-1}
  \varphi(y_j) \right| +
  \\
  + \limsup_{n\to\infty} \frac{1}{n}\sum_{j=k}^{n-1}
  \left| \varphi(y_j) - \varphi(x_j)
  \right| +
  \\
  + \limsup_{n\to\infty} \left| \frac{1}{n}\sum_{j=0}^{n-1}
  \varphi(x_j) - \int_{\mcM} \varphi \, d\mu \right|,
  \end{multline*}
  where $(y_n)$ is the discretization of the path
  $\gamma\in\Gamma_y$, and $(x_n)$ is its $\Phi_{y,x}$-image. We
  know that the measures $F_{*}^{\delta} W_x|_{E_x'}$ and $(F^{\delta}
  \circ \Phi_{y,x})_*
  W_y|_{E_y'}$ are absolutely continuous; thus, the image of a
  typical sequence is a typical sequence. Hence, for a typical
  trajectory $\gamma\in E_y'$ the first and the last summands do not
  exceed $\varepsilon$. So, we have
  \begin{equation}\label{eq:integral_estimate2}
  \limsup_{n\to\infty} \left| \frac{1}{n\delta} \int_0^{n\delta}
  \varphi(\gamma(t)) \, dt - \int_{\mcM} \varphi \, d\mu \right|
  <2\varepsilon
  \end{equation}
  for a typical path $\gamma\in E_y'$, which implies, that
  \begin{equation}\label{eq:integral_estimate3}
  \limsup_{T\to\infty} \left| \frac{1}{T} \int_0^{T}
  \varphi(\gamma(t)) \, dt - \int_{\mcM} \varphi \, d\mu \right|
  <2\varepsilon.
  \end{equation}

  For a typical path from $E_y'$ we have obtained an
  estimate on the difference between the integral of $\varphi$ and its
  average along the path. Let us extend this statement to all the
  $E_{y}$. Namely, repeating the arguments used in the
  proof of Remark~\ref{r:typical}, we see that almost every path
  from $E_{y}$ can be decomposed into a finite starting
  segment and a path from some~$E_{z}'$ for some $z$
  close to $M$. Then, the estimate~\eqref{eq:integral_estimate3} holds
  also for a $W_y$-typical path from $E_y$.

  But the definition of $E_y$ does not depend on~$\delta$,
  thus, choosing arbitrarily small $\delta$ and $r$, we have
  finally:
  $$
  \limsup_{n\to\infty} \left| \frac{1}{T} \int_0^{T}
  \varphi(\gamma(t)) \, dt - \int_{\mcM} \varphi \, d\mu \right|=0.
  $$
  Hence,
  $$
  \frac{1}{T}\int_0^T \varphi(\gamma(t))\, dt \xrightarrow[T\to\infty]{}
  \int \varphi \, d\mu
  $$
  for a typical path $\gamma\in E_{y}$.

  Recall that by definition the measures $\mu_t$ weakly
  converge to $\mu$ if and
  only if for every continuous function $\varphi$ we have
  $$
  \int \varphi \, d\mu_t \to \int \varphi \, d\mu.
  $$
  Moreover, it suffices to check such convergence for a well-chosen
  countable family~$\varphi_k$. We have already obtained, that for
  any function $\varphi$ and for a typical path $\gamma\in E_{y}$
  $$
  \left(\frac{1}{T} \gamma_* \leb_{[0,T]}\right) (\varphi)=
  \frac{1}{T}\int_0^T \varphi(\gamma(t))\, dt \xrightarrow[T\to\infty]{}
  \int \varphi \, d\mu.
  $$

  A countable family of typically satisfied conditions still is a
  typically satisfied condition, and hence for $W_y$-almost
  every trajectory $\gamma\in E_y$
  $$
  \lim\limits_{t\to\infty} \frac{1}{t} \gamma_*\leb_{[0,t]} =\mu.
  $$
  This completes the proof of the lemma.
\end{pfof}


\subsection{Codimension higher than one}\label{ss:higher_codimension}

Here we present a construction which permits us to handle the case
of transversely conformal foliation of codimension higher than one.
Let $\mcF$ be such a foliation. Equip $M$ with a Riemannian metric
and for any point $x\in M$, let a transversal $\mcG_x$ be an image
of the image of the exponential map for a small disk in
$(T_x\mcF)^{\perp}$ disk. We notice that these transversals depend
smoothly on $x$, and in a small neighborhood of every point $x$ for
every point $y$ there exists a unique point $z$ in its small
$\mcF$-leafwise neighborhood, such that $z\in\mcG_y$.

We remark that the family of transversals $\{\mcG_x\}_{x\in M}$ does
not necessarily form a foliation, for the reason that we can not
control intersections of transversals with starting points on
different $\mcF$-leaves.

Now, for a point $x_0\in\mcM$, let us consider the set
$$
\overline{M}=\cup_{y\in \widetilde{\mcF}_{x_0}} \{y\}\times \mcG_y.
$$
Note that $\overline{M}$ is a manifold with boundary, naturally
inheriting from $M$ its foliation structure (except for the fact
that for $\overline{M}$ some leaves intersect the boundary $\partial
\overline{M}$) and Riemannian metric on the leaves. But, on
$\overline{M}$ we have a natural smooth transversal foliation
$\mcG$, leaves of which are $\{y\}\times \mcG_y$. Now, we can apply
the same arguments as these used in the codimension one case, to
prove Lemma~\ref{l:positive}.

\subsection{Non-positivity of Lyapunov exponents}

This section is devoted to the following lemma:
\begin{lemma}\label{l:non-positive}
Let $\mcF$ be a transversely conformal foliation, $\mcM$ be a
minimal set in~$\mcF$, and $\mu$ be a harmonic ergodic measure
supported on~$\mcM$. Then, $\lambda(\mu)\le 0$.
\end{lemma}
\begin{pf}
We will present the proof for the case of codimension one foliation,
using the existence of a transversal foliation~$\mcG$. Then, it is
generalized to an arbitrary codimension case in the same way as in
the proof of Lemma~\ref{l:positive}.

Assume the contrary: let $\lambda(\mu)>0$. We take some
$\alpha,\beta$, $\alpha<\beta<\lambda(\mu)$. Prove first that the
measure $\mu$ does not charge any leaf. Assume the contrary:
$\mu(L)>0$ for some leaf $L$. Then, the density with respect to the
volume $\frac{d\mu|_{L}}{d\vol_g}$ is a harmonic function on the
leaf~$L$. Moreover, this function is positive, bounded (because of
harmonicity and boundedness of geometry of~$L$) and of integral~$1$.
Extending this function by 0 to the complementary, we obtain a
harmonic measurable positive leafwise integrable function.
Garnett~\cite[Proposition~1, p.~295]{Garnett} have proved that such
a function should be leafwise constant. Thus, it is equal to a
positive constant on~$L$ and hence (as it is integrable) $L$ is a
compact leaf. But for a compact leaf the Lyapunov exponent
equals~$0$, as the corresponding Dirac measure is transversely
invariant.

Let us choose a point $x_0\in\mcM$, typical in the sense of Lyapunov
exponents: for $W_{x_0}$-almost every path $\gamma\in\Gamma_{x_0}$
the corresponding path has the Lyapunov exponent equal
to~$\lambda(\mu)$.

First, let us consider the simplest case: suppose, that the
foliation $\mcG$ preserves the metric~$g$. Also, we suppose that any
$\mcF$-along holonomy extends to some fixed neighborhood $U\subset
\mcG_{x_0}$. Finally, we suppose that the measure $\mu$ does not
charge any leaf, or equivalently, that measures $\nu_{\cdot}$ have
no points of positive measure.

In this case, the measure $\mu$ induces on every leaf $\mcG_x$ a
conditional measure $\nu_x$, which is harmonic in the sense of
measure-valued functions (see Section~\ref{ss:non-divergence}). Let
us take some $T>0$ and $C>1$ and for every path $\gamma$, starting
at~$x_0$, try to find a $\tau=\tau(\gamma)$, $T<\tau<2T$, as a
minimal value $t_0$ in this interval possessing the following
property:
\begin{equation}\label{eq:expansion}
\forall t\in[0,\tau] \quad
h'_{\gamma|_{[\tau,t]}}(\gamma(\tau))<Ce^{-\beta(\tau-t)}.
\end{equation}
Here, we use $h_{\gamma|_{[\tau,t]}}$ as a short notation for
$h_{\gamma|_{[t,\tau]}}^{-1}$ ($t<\tau$, and thus the first notation
is not absolutely clear). Note, that as the holonomy is taken in the
inverse sense, from the moment $\tau$ to $t<\tau$, one can expect
that the derivative is will be small.

Note that the non-existence of such $\tau$ means, that
$$
\frac{h'_{\gamma|_{[0,2T]}}}{h'_{\gamma|_{[0,T]}}} \le e^{\beta T}
(x_0),
$$
so for $W_{x_0}$-almost every path $\gamma$ and for every $T$,
sufficiently big such $\tau$ exists. If in the interval $[T,2T]$ we
can not find a moment satisfying~\eqref{eq:expansion}, then we
choose $\tau(\gamma)=2T$. Finally, we remark that $\tau(\cdot)$ is a
Markovian moment.

The transversal measures $\nu_x$ depend on $x$ in a harmonic way,
thus for a transversal interval $I\subset \mcG_{x_0}$ the measures
of its holonomy images $\nu_{\gamma(t)}(h_{\gamma|_{[0,t]}}(I))$
form a martingale. The expectation of value of this martingale at
the Markovian moment~$\tau$ should be equal to its initial value:
\begin{equation}\label{eq:transversal}
\Expect \nu_{\gamma(\tau(\gamma))}
(h_{\gamma|_{[0,\tau(\gamma)]}}(I)) = \nu_{x_0}(I).
\end{equation}
Now, note that due to the definition of $\tau(I)$ for all the paths
$\gamma$ with $\tau(\gamma)<2T$ a neighborhood $V$ of $\gamma(\tau)$
is contracted exponentially by the holonomy
$h_{\gamma|_{[\tau,0]}}$, the radius of $V$ is bounded from below by
means of $C$, $\alpha, \beta$ and the geometry of the foliation.
Namely,
$$
|h_{\gamma|_{[\tau,0]}}(U_{\varepsilon}(\gamma(\tau)))|\le
e^{-\alpha \tau} \le e^{-\alpha T},
$$
where $\varepsilon>0$ does not depend on~$T$. Passing from inverse
to direct time we see that an exponentially small neighborhood of
$x_0=\gamma(0)$ is expanded:
$$
h_{\gamma|_{[0,\tau]}}(U_{e^{-\alpha T}}(x_0)) \supset
U_{\varepsilon}(\gamma(\tau)).
$$
Now, let us estimate the left part of~\eqref{eq:transversal} for
$I=U_{e^{-\alpha T}}(x_0)$:
\begin{multline}\label{eq:growth}
\Expect \nu_{\gamma(\tau(\gamma))}
(h_{\gamma|_{[0,\tau(\gamma)]}}(I)) \ge
\\
\ge \int_{\{\gamma: \tau(\gamma)<2T\}} \nu_{\gamma(\tau(\gamma))}
(h_{\gamma|_{[0,\tau(\gamma)]}}(I)) \, dW_{x_0}\ge
\\
\ge \int_{\{\gamma: \tau(\gamma)<2T\}}
\nu_{\gamma(\tau(\gamma))}(U_{\varepsilon}(\gamma(\tau))) \,
dW_{x_0} \ge
\\
\ge W_{x_0} \{\gamma: \tau(\gamma)<2T\} \cdot
\inf\limits_{y\in\mcM} (\nu_y(U_{\varepsilon}(y))).
\end{multline}
As $\mcM$ is compact and $\supp\mu=\mcM$, the infimum in the last
term of~\eqref{eq:growth} is positive. Thus, the last term stays
separated from 0 as $T$ tends to infinity, so it is no less than
some constant $c_0>0$. Hence~\eqref{eq:transversal}
and~\eqref{eq:growth} imply that
$$
\nu_{x_0}(U_{e^{-\alpha T}}(x_0)) \ge W_{x_0} \{\gamma:
\tau(\gamma)<2T\} \cdot \inf\limits_{y\in\mcM}
(\nu_y(U_{\varepsilon}(y)))>c_0,
$$
and thus the left term does not tend to~$0$ as $T$ tends to
infinity. This contradicts the fact that the measure~$\nu_{x_0}$ can
not have atoms. We have obtained the desired contradiction. So, this
case is handled.

Let us now consider the case of generic Riemannian structure (not
necessarily preserved by the transversal foliation). Note, that the
harmonic measure $\mu$ still defines conditional measures on the
transversals $\{\mcG_y\}$, which are its Fubini conditional measures
with respect to $\vol_g$ on the leaves. We still suppose that the
$\mcF$-along holonomy maps are defined on the entire transverse
leaves $\{\mcG_{y}\}$. Also, we add the following (simplifying the
explanation of this step) hypothesis: all the leaves of $\mcF$ are
simply connected.

For these conditional distributions, the harmonicity condition
implies that for a transversal interval $I$ at a point $x\in M$ and
a function $\rho$, we have:
\begin{equation}\label{eq:integral}
\int_{I} \rho(y) \, d\nu_x (y) =\Expect
\int_{h_{\gamma|_{[0,t]}}(I)}
\frac{p(h_{\gamma|_{[t,0]}}(z),z;t)}{p(x,\gamma(t);t)}
\rho(h_{\gamma|_{[t,0]}}(z)) \, d\nu_{\gamma(t)}(z),
\end{equation}
where the expectation is taken in the sense of $W_x$. To prove this
formula, we take a smooth function $f$ on $M$ supported in the
neighborhood of $I$, with integral on $\mcF_y$ equal to $\rho(y)$.
For every fixed $t$, as the support of $f$ tends to $I$, the
integral of $f$ with respect to $\mu$ tends to the left hand side
of~\eqref{eq:integral}, and the integral of $D^t f$ to the right
hand side. The harmonicity of $\mu$ implies that these two integrals
coincide, which proves the formula.

Now, let us repeat the arguments used in the similar case with the
following modification: we consider only discrete moments of time
$t=k\delta$, where sufficiently small $\delta>0$ is fixed.

Applying~\eqref{eq:integral} several times for the initial function
$\rho=\mathbf{1}_{I}$, we obtain, that for a Markovian moment
(taking discrete values) $t(\gamma)=k(\gamma)\delta$
\begin{equation}\label{eq:expectation}
\nu_x (I) =\Expect \int_{h_{\gamma|_{[0,t]}}(I)}
\prod_{j=1}^{k}\frac{p(z_{j-1},z_j;\delta)}{p(x_{j-1},x_j;\delta)}
\, d\nu_{\gamma(t)}(z),
\end{equation}
where $x_j=\gamma(j\delta)$,
$z_j=h_{\gamma|_{[k\delta,j\delta]}}(z)$, and the expectation is
taken in the sense of the measure $W_{x_0}$.

Let us take, as in the previous case, for $T=K\delta$ sufficiently
big, a Markovian moment $\tau(\gamma)=k(\gamma)\delta$ defined as
the smallest value in the interval $[T,2T]$ such that for every
$t=l\delta<\tau$
$$
h'_{\gamma|_{[\tau,t]}}(\gamma(\tau))<Ce^{-\beta(\tau-t)}.
$$
Once again, if such a moment does not exist, we take $\tau=2T$.

For $T$ sufficiently big, the probability of $\tau<2T$ is close
to~1. Now, let us take $I=U_{e^{-\alpha T}}(x_0)$. For most paths,
as we know, the holonomy maps expand exponentially and thus
$h_{\gamma|_{[0,\tau]}}(I)\supset U_{\varepsilon}(\gamma(\tau))$.

Note that for most of the paths starting at $x_0$, and for a point
$z_0$ in the holonomy preimage
$h_{\gamma|_{[\tau,0]}}(U_{\varepsilon}(\gamma(\tau)))$, the product
of the quotients of the heat kernels (due to estimates analogue to
these of Lemma~\ref{l:positive}) is bounded from below by some
constant $c_1>0$. Let us denote
\begin{multline*}
N=\Bigl\{ \gamma\in\Gamma_{x_0}: \tau(\gamma)<2T,
h_{\gamma|_{[0,\tau]}}(I)\supset U_{\varepsilon}(\gamma(\tau)),
\\
\forall z_0\in h_{\gamma|_{[\tau,0]}}(U_{\varepsilon}(\gamma(\tau)))
\quad \prod_{j=1}^{k}\frac{p(z_{j-1},z_j;\delta)}
{p(x_{j-1},x_j;\delta)} \ge c_1
 \Bigr\}
\end{multline*}

Thus, the right hand side of~\eqref{eq:expectation} can be estimated
as
\begin{multline}\label{eq:expansion_estimate}
\Expect \int_{h_{\gamma|_{[0,t]}}(I)}
\prod_{j=1}^{k}\frac{p(z_{j-1},z_j;\delta)}{p(x_{j-1},x_j;\delta)}
\, d\nu_{\gamma(t)}(z)\ge
\\
\ge \int_N \int_{h_{\gamma|_{[0,t]}}(I)}
\prod_{j=1}^{k}\frac{p(z_{j-1},z_j;\delta)}{p(x_{j-1},x_j;\delta)}
\, d\nu_{\gamma(t)}(z) \, dW_{x_0} \ge
\\
\ge \int_N \int_{U_{\varepsilon}(\gamma(\tau))} c_1 \,
d\nu_{\gamma(t)}(z) \, dW_{x_0} =
\\
= \int_N c_1 \nu_{\gamma(t)}(U_{\varepsilon}(\gamma(\tau))) \,
dW_{x_0} \ge W_{x_0}(N) \cdot c_1 \cdot c_0.
\end{multline}
Once again we see that the measure $\nu_{x_0}(I)$ does not tend to 0
as $I$ contracts to $x_0$, which contradicts the fact that the
measure $\nu$ can not have atoms.

The higher codimension case is handled in the same way by working in
$\overline{M}$ defined in Section~\ref{ss:higher_codimension}. The
measure $\mu$ defines a $\sigma$-finite measure on $\overline{M}$,
which harmonic in the sense of integral definition for compactly
supported test function. Considering leafwise Brownian motion in
$\overline{M}$ (with the possibility of exiting through the
boundary) we see that this measure is superharmonic in the sense
that $\mu\ge D^t_* \mu$, and the same is true for the conditional
measures $\nu_x$.

Note that as in~\eqref{eq:expansion_estimate} the only trajectories
used for estimates are those who arrive in the
$\varepsilon$-neighborhood of the $\gamma(\tau)$; hence due to
exponential contraction in the inverse time they stay closer and
closer to the main leaf $\mcF_{x_0}$. In particular, they stay in
$\overline{M}$ for all the time in the interval $[0,\tau]$.

Adding all this together, we see that the same estimates work in
this case.
So, the general case is handled in the same way.
\end{pf}

\section{Acknowledgements}

The authors would like to thank \'E.~Ghys, who introduced us to the problem, 
and F.~Ledrappier for a careful reading and for many interesting comments on this text. 
We also thank S. Crovisier, who explained to us
the Contraction Lemma, and S.~Frankel, A.~Gorodetski, Yu.~Ilyashenko and A.~Navas for fruitful
discussions.


\begin{thebibliography}{Dillo 83}

\bibitem [An] {Antonov} {\sc V. Antonov.} Modeling of processes of cyclic evolution type. Synchronisation of a random signal.
Vestnik Leningrad. Univ. Mat. Mekh. Astronom. (1984), p. 67-76. 

\bibitem [A-K]{Arnold-Krylov} {\sc V. I. Arnold \& A. L. Krylov.} Uniform distribution of points
on a sphere and some ergodic properties of solutions of linear
ordinary differential equations in a complex region. {\em Sov. Math.
Dokl.}, {\bf 4} (1963), pp.~1--5.

\bibitem [As] {Asuke} {\sc T. Asuke.} On transversely flat conformal foliation with good measures II. {\em Hiroshima Math. J.} {\bf 28}
(1998), pp.~523--525.

\bibitem [Ba] {Baxendale} {\sc P. H. Baxendale.} Lyapunov exponents and relative entropy for a stochastic flow of diffeomorphisms.
{\em Probab. Theory Related Fields.} {\bf 81} (1989), pp.~521--554.

\bibitem [C-D] {Calegari} {\sc D. Calegari \& N. M. Dunfield.} Laminations and groups of homeomorphisms of the circle.
{\em Invent. Math.} {\bf 152} (2003), no. 1, pp.~149--204.

\bibitem [Can] {Candel} {\sc A. Candel.} The harmonic measures of L.~Garnett. {\em Adv. Math.}
{\bf 176} (2003) no. 2, pp.~187--247.

\bibitem [C-C] {Candel1} {\sc A. Candel, L. Conlon.} Foliations II,
{\em Graduate Studies in Mathematics}, v. 60, 2003.

\bibitem [Cha] {Chavel} {\sc I. Chavel.} Eigenvalues in Riemannian
geometry, {\em Pure and Applied Mathematics} 115, Academic Press,
Orlando, 1984.

\bibitem [Ca] {Carriere} {\sc Y. Carri\`ere.} Flots riemanniens.
{\em Ast\'erisque} no. 116 (1984), pp.~31-52.

\bibitem [C-U] {Clozel-Ullmo}{\sc L. Clozel \& E. Ullmo.} Correspondances modulaires et mesures invariantes. {\em J. reine angew. Math.}
{\bf 558} (2003), p. 47-83.

\bibitem [C-L-Y] {CLY} {\sc S. Y. Cheng, P. Li, S. T. Yau}, On the upper
estimate of the heat kernel of a complete Riemannian manifold. Amer.
J. Math. 103 (1981), no. 5, 1021--1063.

\bibitem [Den] {Denjoy} {\sc A. Denjoy.} Sur les courbes d\'efinies par des
\'equations diff\'erentielles \`a la surface du tore. \emph{J. Math.
Pure et Appl.}, \textbf{11} (1932), pp.~333--375.

\bibitem [De] {Deroin} {\sc B. Deroin.} Hypersurfaces Levi-plates immerg\'ees dans les surfaces
complexes de courbure positive. {\em Ann. Scient. \'Ec. Norm. Sup.}
(2005).

\bibitem [D-K-N] {D-K-N} {\sc B. Deroin \& V. Kleptsyn \& A. Navas.} Sur la dynamique unidimensionnelle 
en r\'egularit\'e interm\'ediaire. Preprint 2005. 

\bibitem [Ep] {Epstein} {\sc D. B. Epstein.} Transversely hyperbolic
$1$-dimensional foliations. {\em Ast\'erisque} no. 116 (1984),
pp.~53--69.

\bibitem [Fe] {Fenley} {\sc S. Fenley.} Foliations, topology and geometry of $3$-manifolds:~${\bbR}$-covered foliations
and transverse pseudo-Anosov flow. {\em Comment. Math. Helv.} {\bf 77} (2002), no.3, p.~415--490.

\bibitem [F-G] {Feres-Ghani} {\sc R. Feres \& A. Zeghib.} Dynamics of
the space of harmonic functions and foliated Liouville problem. To
appear in {\em Ergod. Theory Dyn. Syst.}.

\bibitem [F-S] {Fornaess-Sibony} {\sc J. E. Fornaess \& N. Sibony.} Harmonic currents of finite energy and laminations.
arXiv:math.CV/0402432v1.

\bibitem [F] {Furman} {\sc A. Furman.}
Random walks on groups and random transformations. Handbook of
dynamical systems, Vol. 1A, 931--1014, North-Holland, Amsterdam,
2002.

\bibitem [Fu1]{Furstenberg} {\sc H. Furstenberg.} Non commuting random matrices product. {\em Trans. Amer. Math. Soc.}
(1963) {\bf 108}, p. 377-428.

\bibitem [Fu2]{Furstenberg-minimality} {\sc H. Furstenberg.} Strict
ergodicity and transformation of the torus. {\em Amer. J. Math.}
{\bf 83} (1961) pp.~573--601.

\bibitem [F-K] {Furstenberg-Kesten} {\sc H. Furstenberg \& H. Kesten.}
Product of random matrices. {\em Ann. Math. Stat.}, {\bf 31} (1960),
pp.~457--469.

\bibitem [Ga]{Garnett} {\sc L. Garnett.} Foliations, the ergodic theorem and Brownian motion.
{\em J. Funct. Anal.} {\bf 51}  (1983), no. 3, p. 285-311.

\bibitem [Gh1] {Ghys1}
 {\sc \'E. Ghys.}  Flots transversalement affines et tissus feuillet\'es.
{\em Mém. Soc. Math. France} No. 46 (1991), p. 123-150.

\bibitem [Gh2] {Ghys2}
 {\sc \'E. Ghys.}  Sur l'uniformisation des laminations paraboliques.
 In: Integrable systems and foliations/Feuilletages et syst\`emes int\'egrables
 (Monpellier, 1995), pp.~73--91.

\bibitem [Gh3] {Ghys3}
 {\sc \'E. Ghys.} Gauss--Bonnet theorem for 2-dimensional foliations,
 J. Funct. Anal., {\bf 77} (1988), no.~1, pp.~51--59.

\bibitem [Gh4] {Ghys4} {\sc \'E. Ghys.} Laminations par surfaces de Riemann.
\emph{Dynamique et g\'eom\'etrie complexes (Lyon 1997)}, ix, xi,
pp.~49--95, Panor. Synth\`eses, 8, \emph{Soc. Math. France, Paris},
1999.

\bibitem [Gui] {Guivarc'h} {\sc Y. Guivarc'h.} Quelques
propri\'et\'es asymptotiques pes produits de matrices al\'eatoires.
\'Ecole d'\'et\'e (Saint Flour, 1978), pp.~177-250.

\bibitem [G-R]{Guivarc'h-Raugi} {\sc Y. Guivarc'h \& A. Raugi.}
Sur les mesures invariantes de certaines cha\^{\i}nes de Markov
d\'efinies par des transformations homographiques. Random walks and
stochastic processes on Lie groups, p. 62-65, Inst. \'Elie Cartan,
Nancy, 1983.

\bibitem [Hae]{Haefliger} {\sc A. Haefliger} Stuctures feuillet\'ees
et cohomologie \`a valeurs dans un faisceau de groupo\"ides,
\emph{Comment. Math. Helv.}, \textbf{32} (1958), pp.~248-329.


\bibitem [Ham1]{Hamenstadt1} {\sc U. Hamenst\"adt.} Harmonic measures for compact negatively curved manifolds.
{\em Acta Math.} {\bf 178} (1997), no. 1, pp. 39. 

\bibitem [Ham2] {Hamenstadt2} {\sc U. Hamenst\"adt.} Positive eigenfunctions on the universal covering of a compact 
negatively curved manifold. Preprint. 

\bibitem [He] {Helmberg} {\sc G. Helmberg.}
A theorem on equidistribution on compact groups. {\em Pacific J.
Math.} {\bf 8} (1958), p. 227-241.

\bibitem [Hur] {Hurder} {\sc S. Hurder.} Exceptional minimal sets
for $C^{1+\alpha}$-group actions on the circle. \emph{Ergodic Theory
and Dynamical Systems}, \textbf{11} (1991), pp.~455--467.

\bibitem [K-M] {K-M} {\sc V. A. Kaimanovich, H.~Mazur.} The Poisson boundary of the
mapping class group, {\em Invent. Math}, {\bf 125} (1996),
pp.~221--264.

\bibitem [Kai1] {Kaimanovich} {\sc V. A. Kaimanovich.}  Brownian motion on foliations: entropy,
invariant measures, mixing, Funct. Anal. Appl., {\bf 22} (1988), no.
4, pp. 326-328.

\bibitem [Kai2]{Kaimanovich2} {\sc V. A. Kaimanovich.} Brownian
motion and harmonic functions on covering manifolds. An entropy
approach. {\em Dokl. Akad. Nauk. SSSR} {\bf 288} (1986),
pp.~1045--1049 (Russian); English translation: Soviet Math. Dokl.,
{\bf 33} (1986), pp.~812-816.

\bibitem [Kaij]{Kaijser} {\sc T. Kaijser.} On stochastic perturbations of iterations of circle maps.  
{\em Phys. D} {\bf 68}  (1993),  no. 2, p. 201-231.


\bibitem [Kak] {Kakutani} {\sc S. Kakutani.} Random ergodic theorems and Markoff processes with a stable distribution.
{\em Proceedings of the Second Berkeley Symposium on Mathematical Statistics and Probability, 1950, University of California Press,
Berkeley and Los Angeles} (1951), pp.~247--261.

\bibitem [K-H] {Katok-Hasselblatt} {\sc A. Katok \& B. Hasselblatt.}
Introduction to the modern theory of dynamical systems. {\em
Cambridge University Press, Cambridge} (1995).

\bibitem [Kea] {Keane}
{\sc M. Keane.} Non-ergodic interval exchange transformations.  {\em
Israel J. Math.} {\bf 26} (1977), no. 2, p. 188-196.

\bibitem [Ke-N] {Keynes-Newton}  {\sc H. B. Keynes \& D. Newton.}
A ``minimal'', non-uniquely ergodic interval exchange
transformation. {\em Math. Z.}  {\bf 148}  (1976), no. 2, p.
101-105.

\bibitem [Kl-N]{Kleptsyn-Nalski} {\sc V. Kleptsyn \& M. Nalski.} Contraction of
orbits in random dynamical systems on the circle. {\em Funct. Anal.
Appl.}, {\bf 38} (2004), no. 4, pp. 36--54.

\bibitem [Led1] {Ledrappier} {\sc F. Ledrappier.} Ergodic properties of the stable foliations.
{\em Ergodic Theory and Related Topics III (1990,G\"ustrow)}, pp.~131-145,
Lecture Notes in Math., 1514.

\bibitem [Led2] {Ledrappier2} {\sc F. Ledrappier.} Application of dynamics to compact manifolds of negative curvature.
{\em Proc. of the International Congress of Mathematicians.} {\bf 1}, no 2, (Z\"urich 1994) pp. 1195-1202, Birkh\"auser, Basel, 1995.

\bibitem [LJ] {Le Jan}
 {\sc Y. Le Jan.} \'Equilibre statistique pour les produits de diff\'eomorphismes al\'eatoires ind\'ependants.
{\em Ann. Inst. Henri Poincar\'e.} {\bf 23} (1987), no. 1,
pp.~111--120.

\bibitem [Li] {Lin} {\sc M. Lin.} On the ``zero-two'' law for
conservative Markov processes. {\em Z. Wahrsch. Verw. Gebiete} {\bf
61} (1982), no. 4, pp.~513--525.

\bibitem [Ma] {Malliavin} {\sc P. Malliavin.} Diffusion et G\'eom\'etrie Diff\'erentielle Globale. Centro Internationale Matematico Estivo.
Varenne, France, ao\^ut 1975.

\bibitem [Pe] {Petite} {\sc S. Petite.} On invariant measures of finite affine type tilings. Preprint.

\bibitem [Pl]{Plante} {\sc J. Plante.} Foliations with measure preserving holonomy.
{\em Ann. of Math.} {\bf 102} (1975), pp.~327--361.

\bibitem [P-G] {PG} {\sc J. F. Plante, S. E. Goodman.} Holonomy and
averaging in foliated sets, J. Differential Geometry, {\bf 14}
(1979), no. 3, pp. 401--407.

\bibitem [R-S]{Ruelle-Sullivan} {\sc D.Ruelle \& D. Sullivan.} Currents, flows, and diffeomorphisms. {\em Topology} {\bf 14}
(1975), pp.~319--327.

\bibitem [Sa] {Sacksteder} {\sc R. Sacksteder.} Foliations and pseudogroups. {\em Amer. J. Math.} {\bf 87} (1965), p.
79-102.

\bibitem [Sch]{Schwartz} {\sc A. Schwartz.} A generalization of
Poincar\'e-Bendixon theorem to closed two dimensional manifolds.
\emph{Amer. J. Math.} \textbf{85} (1963), 453--458.

\bibitem [Scm]{Schwartzman}{\sc S. Schwartzman.} Asymptotic cycles. {\em Ann. of Math.} {\bf 66} (1957), pp.~270--284.

\bibitem [Su1] {cycle} {\sc D. Sullivan.} Cycles for the dynamical study of foliated manifolds and complex manifolds.
{\em Invent. Math.} {\bf 36}  (1976), p. 225--255.

\bibitem [Su2] {Su2} {\sc D. Sullivan.} Conformal dymanical systems.
{\em Lecture Notes in Mathematics}, {\bf 1007} (1983), Springer, New
York, pp. 725--752.

\bibitem [Su3]{Su3} {\sc D. Sullivan.}
Linking the universalities of Milnor-Thurston, Feigenbaum and
Ahlfors-Bers. Topological methods in modern mathematics, (Stony
Brook NY, 1991), p. 543-564.


\bibitem [Ta] {Tarquini} {\sc C. Tarquini.} Feuilletages conformes,
{\em Ann. Inst. Fourier} {\bf 54} (2004), no. 2, pp.~453--480.

\bibitem [Th] {Thurston} {\sc W. Thurston.} Three-manifolds, Foliations and Circles,
II. Unfinished manuscript, 1998.

\bibitem [Ve] {Vershik}
 {\sc A. Vershik.} Polymorphisms, Markov processes,
random perturbations of K-automorphisms. math.DS/0409492.

\bibitem [Yue] {Yue} {\sc C.-B. Yue.} Brownian motion on Anosov foliation and manifold of negative curvature. {\em J. Differential Geom.}
{\bf 41} (1995), pp.~159-183.

\end{thebibliography}
\end{document}